\documentclass[12pt]{amsart}
\pdfoutput=1
\usepackage{amsmath}

\ifx\pdfoutput\undefined
 \usepackage{epsfig}
 \else
 \usepackage[pdftex]{epsfig}
 \usepackage{epstopdf}
 \fi

\theoremstyle{definition}

\newcommand{\edash}{\frac{\ \ }{}}

\begin{document}

\parskip 6pt
\parindent 12pt
\baselineskip 14pt

\title[tessellation\dots]{A tessellation for algebraic surfaces in $\mathbb C\mathbf P^3$ }

\author[Hanson]{Andrew J.~Hanson}
\author[Sha]{Ji-Ping Sha}
\address{Department of Computer Science, Indiana University,
Bloomington, IN 47405} \email{hansona@indiana.edu}
\address{Department of Mathematics, Indiana University,
Bloomington, IN 47405} \email{jsha@indiana.edu}

\maketitle

\bigskip

In this paper we present a systematic and explicit algorithm for
tessellating the algebraic surfaces (real 4-manifolds) $F_n$ in
$\mathbb C\mathbf P^3$ defined by the equation
\begin{equation}
z_0^{\,n}+z_1^{\,n}+z_2^{\,n}+z_3^{\,n}=0 \label{eqn:surface}
\end{equation}
in the standard homogeneous coordinates $[z_0,z_1,z_2,z_3]$, where
$n$ is any positive integer. Note that $F_{4}$ in particular is a
$K3$ surface (see, e.g., \cite{gh78}).

The tessellation we present contains a minimal number of vertices:
they are the $n$-th roots of unity in the six standard projective
lines $\mathbb C\mathbf P^1$ in $\mathbb C\mathbf P^3$, and are the
obvious vertices to start a construction of a natural tessellation
for $F_n$. Our tessellation is invariant under the action of the
obvious isomorphism group of $F_n$ induced by permutations and phase
multiplications of the coordinates, and the action is transitive on
the set of 4-cells. The tessellation is built upon a similar
triangulation for the corresponding algebraic curves in $\mathbb
C\mathbf P^2$, and we believe it can be generalized to the
corresponding algebraic hypersurfaces in $\mathbb C\mathbf P^N$ for
$N>3$.

The tessellation is algorithmically programmable: For any given
positive integer $n$, one first lists all the vertices; then all the
edges, faces, 3-cells, and 4-cells can be produced symbolically from
the list of vertices.  One can then, for example, also formulate the
simplicial complex boundary map matrices and compute the homology,
etc., if one wishes.

Explicit representations of geometric objects such as manifolds are
essential for any attempt to create visual images that help expose
their features. While there exist many powerful mathematical methods
that allow the calculation of the geometric and topological
invariants of manifolds, human perception requires the construction
of visual images.  Thus, it can be useful to develop explicit
descriptions of interesting families of manifolds that can be used
in practice to create visual representations and pictures. Such
explicit representations can also in principle be used to clarify
the calculation and understanding of abstract invariants of the
manifolds. Among the classes of geometric objects that have a long
history of interest are the algebraic varieties defined by
homogeneous polynomials in complex projective spaces.  One such
family, the algebraic curves in $\mathbb C\mathbf P^2$ (see, e.g.,
\cite{ajh94}), has recently served the purpose of providing explicit
images of cross-sections of Calabi-Yau spaces, and has been used to
represent the hidden dimensions of string theory \cite{bg99}, for
which very few other methods of producing images are available.
While one might have guessed that the methods used for  $\mathbb
C\mathbf P^2$ could be extended trivially to $\mathbb C\mathbf P^3$
and higher dimensional projective spaces, the problem turns out to
be fairly complex.

Let us now be more precise. We will show the following:

\noindent{\bf Theorem.} {\it For any given positive integer $n$,
there is a tessellation on $F_n$ with $6n^3$ 4-cells. Each 4-cell is
bounded by four pentahedrons. Each pentahedron is a pyramid with one
quadrilateral face and four triangular faces. The tessellation is
invariant under the action of the group $\Gamma_n$, where $\Gamma_n$
consists of isomorphisms of $F_n$ induced from permutations and
phase multiplications of the homogeneous coordinates of $\mathbb
C\mathbf P^3$. The group $\Gamma_n$ acts transitively on the set of
4-cells of the tessellation.}

Altogether, the tessellation has $6n$ vertices, $12n^2$ edges,
$8n^2+7n^3$ 2-cells ($3n^3$ quadrilaterals and $8n^2+4n^3$
triangles), $12n^3$ 3-cells (pyramids) and $6n^3$ 4-cells. It is
known that the Euler characteristic of any smooth algebraic surface
of degree $n$ in $\mathbb C\mathbf P^3$ is $6n-4n^2+n^3$ (see, e.g.,
\cite{gh78}). One handily verifies from our tessellation for $F_n$
that this is equal to $6n-12n^2+(8n^2+7n^3)-12n^3+6n^3$, i.e., the
alternating sum of the numbers of vertices, edges, 2-cells, 3-cells,
and 4-cells.

Notice that the restriction to $F_n$ of the natural projection
$\mathbb C\mathbf P^3 \setminus\{[0,0,0,1]\} \to \mathbb C\mathbf
P^2$, given by $[z_0,z_1,z_2,z_3] \mapsto [z_0,z_1,z_2]$, is a
regular n-fold branched covering
\begin{equation}\label{eqn:branchcover}
\sigma: F_n \to \mathbb C\mathbf P^2 \end{equation} which is
branched over the algebraic curve in $\mathbb C\mathbf P^2$ defined
by the equation
\begin{equation}\label{eqn:curve}
z_0^{\,n}+z_1^{\,n}+z_2^{\,n}=0 \ .
\end{equation}

The tessellation of $F_n$ we present is a lift from $\sigma$ of a
tessellation of $\mathbb C\mathbf P^2$, which is an extension of a
tessellation (triangulation) of the algebraic curve
(\ref{eqn:curve}). This approach greatly reduces the difficulty
caused by the topological complexity of $F_n$, as the geometry and
topology of $\mathbb C\mathbf P^2$ are much easier to handle and
visualize. We also implicitly assume that $\mathbb C\mathbf P^2$ is
equipped with the standard Fubini-Study Riemannian metric. In
particular, every projective line $\mathbb C\mathbf P^1$ in $\mathbb
C\mathbf P^2$ is totally geodesic, and, with the induced metric, is
a round 2-sphere; the real projective planes are also totally
geodesic and have induced metric of constant curvature.

\bigskip

\section{Tessellation of the algebraic curve}
\label{section:tesscurve}

 Denote by $S_n$ the algebraic curve in $\mathbb C\mathbf
P^2$ defined by (\ref{eqn:curve}). In this section, we will
tessellate (i.e., triangulate) $S_n$ in a specific way so that we
can extend the tessellation to the $\mathbb C\mathbf P^2$ in the
next section. The tessellation is in fact a lifting of a natural
tessellation on $\mathbb C\mathbf P^1$ for the given $n$.

The projection $\mathbb C\mathbf P^2\setminus\{[0,0,1]\} \to \mathbb
C\mathbf P^1$, given by $[z_0,z_1,z_2] \mapsto [z_0,z_1]$ induces a
regular $n$-fold branched covering from $S_n$ to the $\mathbb
C\mathbf P^1$ branched at $n$ points,
\begin{equation}\label{eqn:cp1roots}
p_k:=[1,e^{i(\pi+2k\pi)/n}], \ \ \ k=0,\dots,n-1 \ .
\end{equation}

\begin{figure}[!htbp]
\centering
 \centerline{ \psfig{figure=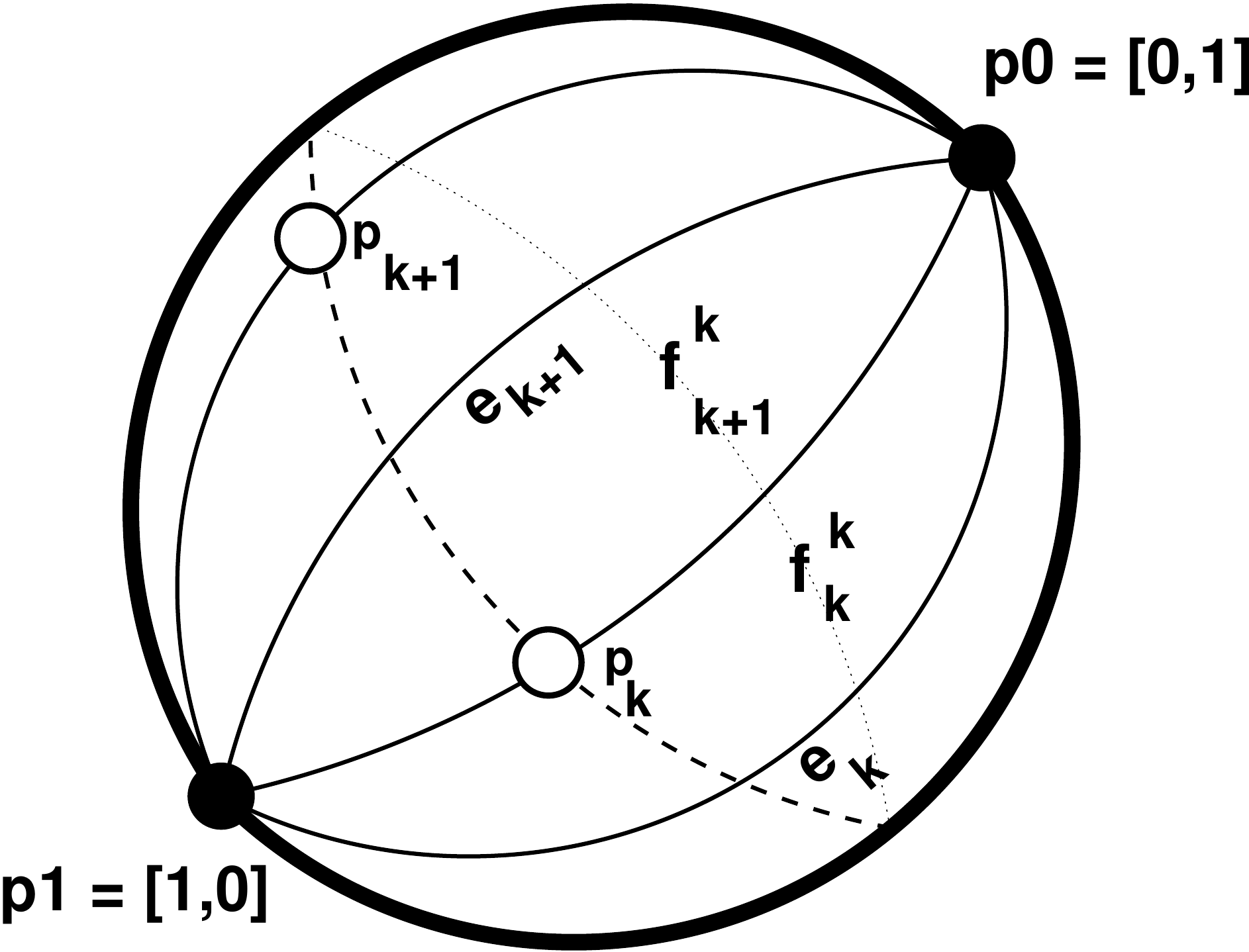,height=4.5cm} }
\caption{}
  \label{cp1.fig}
\end{figure}

We first formulate a tessellation for the $\mathbb C\mathbf P^1$,
which has $n+2$ vertices, $3n$ edges and $2n$ triangles:

Let
\begin{equation}\label{eqn:cp1poles}
p0:=[0,1], \ \ \  p1:=[1,0]
\end{equation}
and join them by the following $n$ paths,
\begin{equation}\label{eqn:p0p1path}
e_k(t):=[\cos t,\sin t e^{i2k\pi/n}], \ 0\le t\le \frac{\pi}{2} \ ,
\ \ k=0,\dots,n-1 \ . \end{equation} Then the $2n$ triangles of the
tessellation are
\begin{equation}\label{eqn:cp1tris}
f_k^k, \ \ f_{k+1}^k \ , \ k=0,\dots,n-1 \ \ \ \text{(mod $n$)} \ ,
\end{equation}
where each of $f_k^k$, $f_{k+1}^k$  is the triangle with vertices
$p0,p1,p_k$, the edge $e_k$, or $e_{k+1}$, respectively, and the
other two edges given by the minimizing geodesics joining $p_k$ and
$p0,p1$ (see Figure \ref{cp1.fig}).

Lifting this triangulation through the branched covering, we then
get a triangulation for $S_n$. There are $3n$ vertices,
$$p0_{k}:=[0,1,e^{i(\pi+2k\pi)/n}], \ p1_{k}:=[e^{i(\pi+2k\pi)/n},0,1], \ p2_{k}:=[1,e^{i(\pi+2k\pi)/n},0],$$
for $k=0,\dots,n-1 \ \ \ \text{(mod $n$)}$, and $2n^2$ triangles. It
is not hard to see that these triangles, as lifts of $f_k^k,
f_{k+1}^k$ and expressed in terms of their vertices, are $\triangle
p0_{j-k}p1_{-j-1}p2_k,  \ \triangle p0_{j-(k+1)}p1_{-j-1}p2_k $,
respectively. We denote them by the following:
\begin{equation}\label{eqn:trilist2}
b_{j-k,-j-1,k} \ , \ \ b_{j-(k+1), -j-1,k} \ , \ \ \ j,k=0,\dots,n-1
\ \text{(mod $n$)} \ .
\end{equation}
To be more clear, we verify the indices in (\ref{eqn:trilist2}) by
showing the edges of these triangles explicitly.

\begin{figure}[!htbp]
\centerline{ \psfig{figure=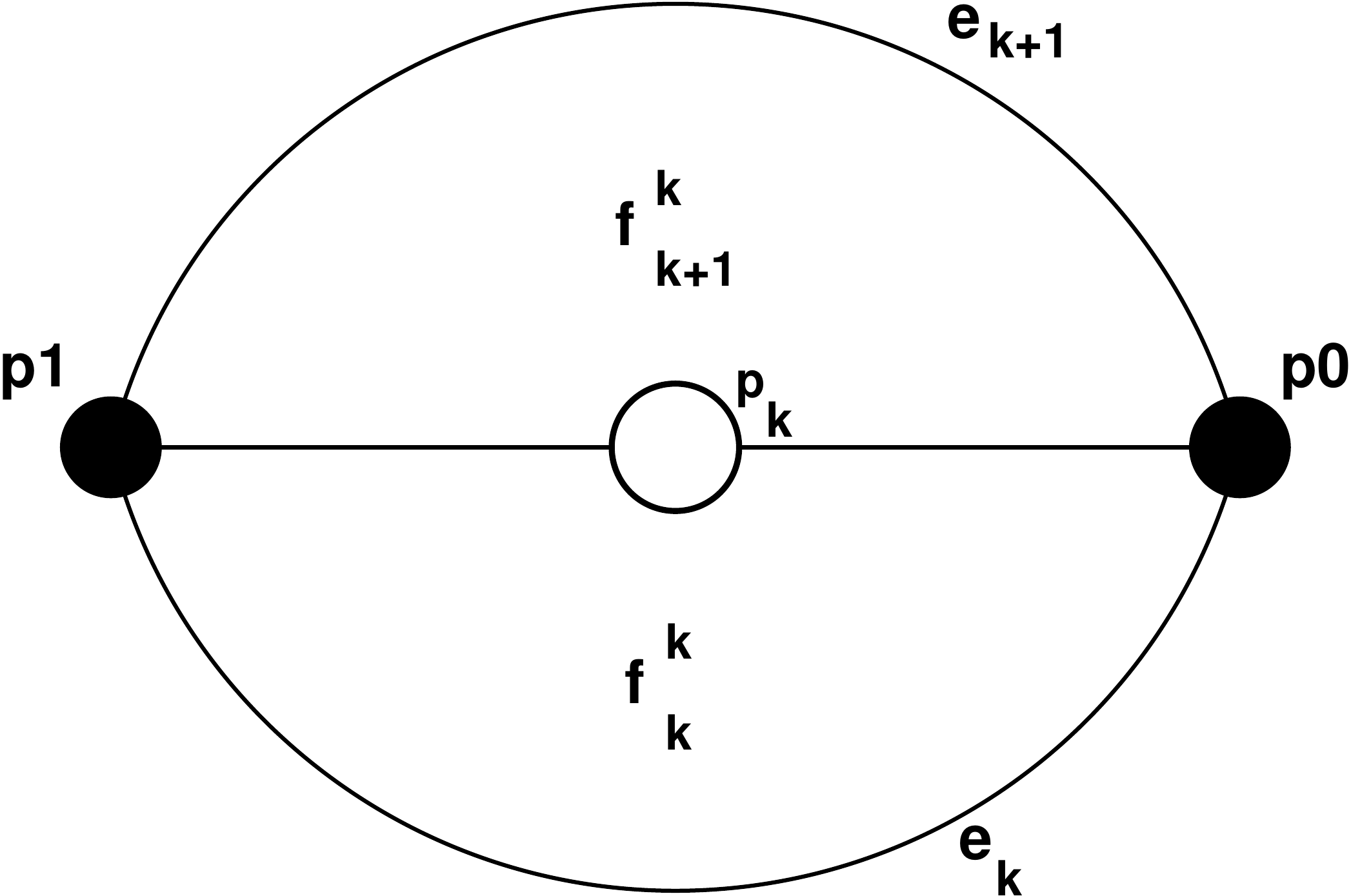,height=4.0cm}
\hspace{.1in}  \psfig{figure=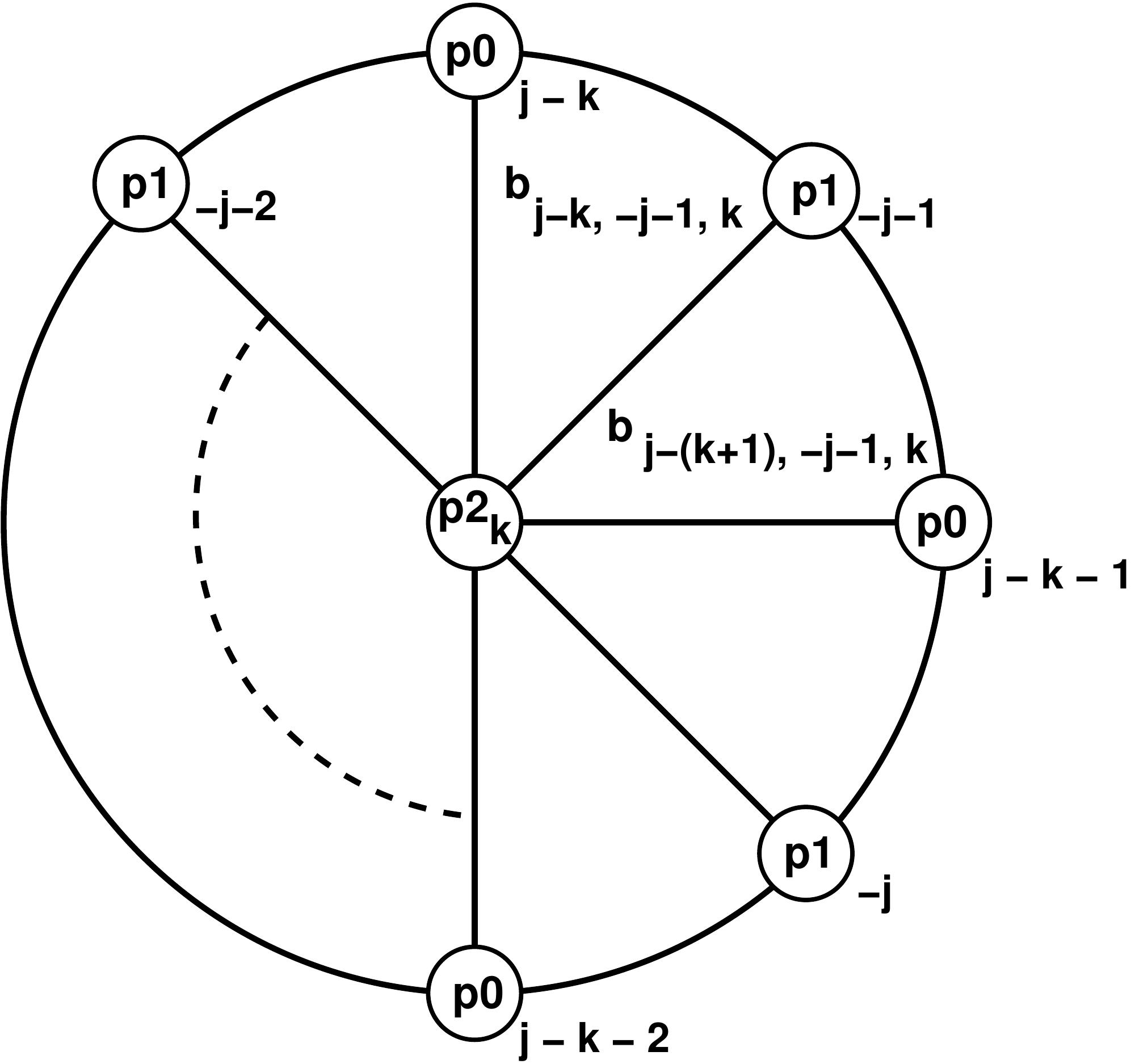,height=6.0cm} }
\caption{}
  \label{branch.fig}
\end{figure}

The three edges of $b_{j-(k+1), -j-1,k}$, in the order $p2_k\edash
p0_{j-(k+1)}\edash p1_{-j-1}\edash p2_k$, can be described as
follows: notice that the first two coordinates give the edges of
$f_{k+1}^k$, in the order $p_k\edash p0\frac{e_{k+1}}{}p1\edash
p_k$, and the factor $e^{i2j\pi/n}$ on the third coordinate
specifies a certain branch to which $f_{k+1}^k$ is lifted.
\begin{eqnarray} \label{eqn:kedge1}
&&[\cos t,\sin t \, e^{i(2k+1)\pi/n},
e^{i2j\pi/n}(-\cos^nt+\sin^nt)^{1/n}] \ , \ \ \pi/4\le t\le\pi/2 \
; \nonumber \\
&&[\sin t,\cos t\, e^{i2(k+1)\pi/n},
e^{i2j\pi/n}(-\sin^nt-\cos^nt)^{1/n}] \ , \ \ 0\le t\le\pi/2
\ ;  \\
&& [\cos t,\sin t \, e^{i(2k+1)\pi/n},
e^{i2j\pi/n}(-\cos^nt+\sin^nt)^{1/n}] \ , \ \ 0\le t\le\pi/4 \ .
\nonumber  \end{eqnarray} Similarly, the three edges of
$b_{j-k,-j-1,k}$, in the order $p2_k\edash p1_{-j-1}\edash
p0_{j-k}\edash p2_k$, as lifts of those of $f_k^k$, in the order
$p_k\edash p1\frac{e_k}{}p0\edash p_k$, are
\begin{eqnarray} \label{eqn:kedge2}
&& [\sin t,\cos t \, e^{i(2k+1)\pi/n},
e^{i2j\pi/n}(-\sin^nt+\cos^nt)^{1/n}] \ , \ \ \pi/4\le t\le\pi/2 \ ;
\nonumber \\
&&[\cos t,\sin t \, e^{i2k\pi/n}, e^{i2j\pi/n}(-\cos^nt-\sin^nt)^{1/n}]
\ , \ \ 0\le t\le\pi/2 \ ;  \\
&&[\sin t,\cos t \, e^{i(2k+1)\pi/n},
e^{i2(j+1)\pi/n}(-\sin^nt+\cos^nt)^{1/n}] \ , \ \ 0\le t\le\pi/4 \ .
\nonumber  \end{eqnarray} Notice that there is a branch shift on the
lift of $p0\edash p_k$  (see Figure \ref{branch.fig}).

From (\ref{eqn:kedge1}) and (\ref{eqn:kedge2}) one then gets the
vertices for the corresponding triangles easily. We note that any
one of the indices for $b$ in (\ref{eqn:trilist2}) is determined by
the other two according to the rule that the sum of the three
indices is equal to $-1$ or $-2$, respectively.

The tessellation is invariant under the isomorphisms of $S_n$
induced from permutations and phase multiplications of the
homogeneous coordinates of $\mathbb C\mathbf P^2$; here a phase
multiplication means multiplying any of the coordinates by a number
of the form $e^{i2k\pi/n}$. To see this, first notice that the
tessellation on $\mathbb C\mathbf P^1$ is obviously invariant under
the corresponding isomorphisms: the vertices are invariant and the
edges are all geodesics while the isomorphisms are isometries. The
tessellation is also obviously invariant under the phase
multiplication of $z_2$ because the latter is just a deck
transformation of the branched covering. Therefore it suffices only
to verify the invariance under interchanging the coordinates $z_1$
and $z_2$.

After interchanging $z_1$ and $z_2$, the three paths in
(\ref{eqn:kedge2}) become
\begin{eqnarray*}
&& [\sin t, e^{i2j\pi/n}(-\sin^nt+\cos^nt)^{1/n},\cos t
e^{i(2k+1)\pi/n}] \ , \ \ \pi/4\le t\le\pi/2 \ ;\\
&&[\cos t, e^{i2j\pi/n}(-\cos^nt-\sin^nt)^{1/n},\sin te^{i2k\pi/n}]
\ , \ \ 0\le t\le\pi/2 \ ; \\
&&[\sin t, e^{i2(j+1)\pi/n}(-\sin^nt+\cos^nt)^{1/n},\cos
te^{i(2k+1)\pi/n}] \ , \ \ 0\le t\le\pi/4 \ .
\end{eqnarray*}
They are the same as
\begin{eqnarray*}
&&[\sin t, e^{i(2j+1)\pi/n}(\sin^nt-\cos^nt)^{1/n},\cos t
e^{i(2k+1)\pi/n}] \ , \ \ \pi/4\le t\le\pi/2 \ ;\\
&&[\cos t, e^{i(2j+1)\pi/n}(\cos^nt+\sin^nt)^{1/n},\sin
te^{i2k\pi/n}] \ , \ \ 0\le t\le\pi/2 \ ; \\
&&[\sin t, e^{i2(j+1)\pi/n}(-\sin^nt+\cos^nt)^{1/n},\cos
te^{i(2k+1)\pi/n}] \ , \ \ 0\le t\le\pi/4 \ ;
\end{eqnarray*}
or
\begin{eqnarray*}
&&[\cos t,
e^{i(2j+1)\pi/n}(\cos^nt+\sin^nt)^{1/n},e^{i2k\pi/n}(\sin^nt)^{1/n}]
\ , \ \ 0\le t\le\pi/2 \ ; \\
&&[\sin t,
e^{i2(j+1)\pi/n}(-\sin^nt+\cos^nt)^{1/n},e^{i2k\pi/n}(-\cos^nt)^{1/n}]
\ , \ \ 0\le t\le\pi/4 \ . \\
&&[\sin t,
e^{i(2j+1)\pi/n}(\sin^nt-\cos^nt)^{1/n},e^{i2k\pi/n}(-\cos^nt)^{1/n}]
\ , \ \ \pi/4\le t\le\pi/2 \ ;
\end{eqnarray*}
which are the edges of $b_{k-(j+1),-k-1,j}$, a lift of $f_{j+1}^j$.
Similarly, interchanging $z_1$ and $z_2$ transforms
$b_{j-(k+1),-j-1,k}$ to $b_{k-j,-k-1,j}$.

It is easy to see that the transformation under these isomorphisms
is transitive on triangles. As the number of the isomorphisms is
$6n^2$, the order of isotropy of each triangle is $3$, consisting of
the cyclic edge permutations. Therefore the transformation is also
transitive on the edges, and obviously on the vertices as well.

We finally point out that the case $n=1$ is somewhat peculiar: the
two triangles share the same three edges. Therefore extra care in
labeling, e.g, specifying the orientation, is needed.

\bigskip

\section{Extended tessellation on $\mathbb C\mathbf
P^2$}\label{section:tesscp2}

In this section, we extend the tessellation of $S_n$ described
in $\S$\ref{section:tesscurve} to a tessellation of the $\mathbb
C\mathbf P^2$. Then, by lifting, that will automatically produce a
tessellation of $F_n$.

Denote the projective line $z_j=0$ by $L_j$, for $j=0,1,2$, and let
\begin{equation}\label{eqn:cp2poles}
p01:=[0,0,1], \ \ p12:=[1,0,0], \ \ p20:=[0,1,0].
\end{equation}

We start by specifying the other $2$-cells for the tessellation.

Note that on $L_0$, the points $p01,p20$ and the intersections with
the $S_n$, namely $p0_k, \ k=0,\dots,n-1$, form the exact same
configuration as (\ref{eqn:cp1poles}) and (\ref{eqn:cp1roots}) on
$\mathbb C\mathbf P^1$ described in $\S$\ref{section:tesscurve}. We
then add the corresponding $2n$ triangles (\ref{eqn:cp1tris});
similarly for the lines $L_1$ and $L_2$. Therefore altogether there
are $6n$ new triangles, which we label as follows:
\begin{equation}\label{eqn:trilist1}
f\!j\,{}_k^k, \ \ f\!j\,{}_{k+1}^k, \ \ \ j=0,1,2 \ \ \text{and} \ \
k=0,\dots,n-1 \ \ \text{(mod $n$)}.
\end{equation}
Label the edges corresponding to those in (\ref{eqn:p0p1path}) by
$ej\,{}_k$. Notice that, for example, as a path, $e1\,{}_k(t)=[\sin
t\, e^{i2k\pi/n},0,\cos t]$.

In the next group, each triangle is formed by minimizing geodesics
joining one of the vertices $p01,p12,p20$  to the edge on $S_n$,
e.g., $p2_j\edash p0_k$ in the case of $p20$. We denote these $3n^2$
triangles as follows:
\begin{equation}\label{eqn:trilist3}
h01_{jk}, \ \ h12_{jk}, \ \ h20_{jk}, \ \ \ j,k=0,\dots,n-1 \ \
\text{(mod $n$)}.
\end{equation}

We remark that all the triangles in (\ref{eqn:trilist3}) are totally
geodesic; one sees, e.g., from (\ref{eqn:kedge1}) that they are
pieces of real projective planes. In fact, all the new $2$-cells we
add will be totally geodesic.

There is one more group of $n^2$ triangles that all have the same
three vertices  $p01,p12,p20$. For clarity, we write down the
following explicit parameterizations for them:
$$g_{jk}(s,t)=[\cos s, \sin s \cos t e^{i2j\pi/n}, \sin s \sin t e^{i2k\pi/n}] \ , \ \ 0\le s,t\le\pi/2 \ ,$$
The three edges of $g_{jk}$ are $e0\,{}_{k-j},e1\,{}_{-k},e2\,{}_j$.
For convenience, we will denote $g_{jk}$ by
\begin{equation}\label{eqn:trilist4}
g_{k-j,-k,j}, \ \ \ \ j,k=0,\dots,n-1 \ \ \text{(mod $n$)}\ ,
\end{equation}
noticing again that any one of the indices of $g$ is determined by
the other two according to the rule that the sum of the three
indices is equal to $0$.

\begin{figure}[!htbp]
 \centerline{
\psfig{figure=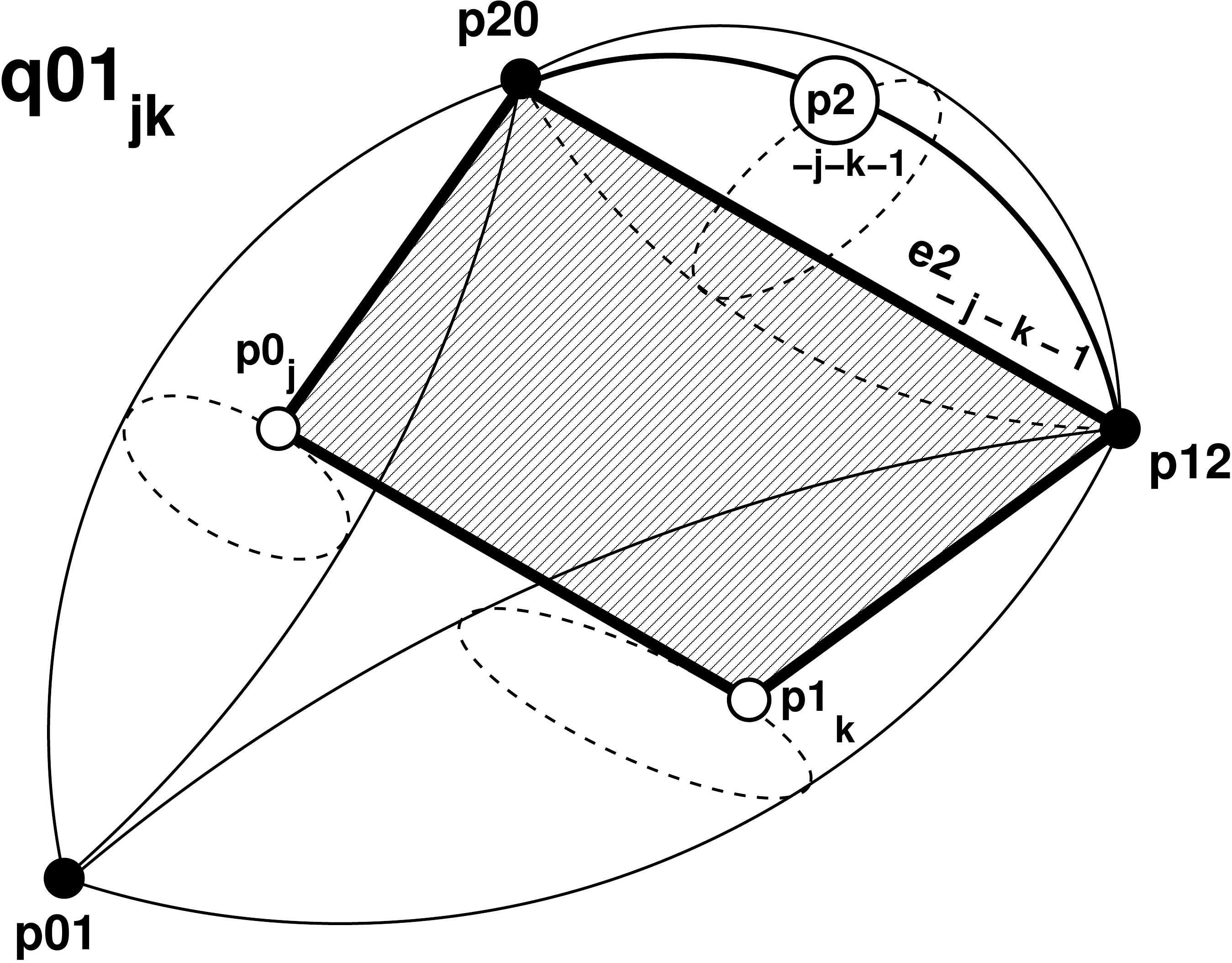,height=5.5cm} }
\caption{}
  \label{quadpic.fig}
\end{figure}

The next set of $2$-cells is a set of $3n^2$ quadrilaterals. They
are in one-to-one correspondence with the edges in $S_n$; each edge
is one side of exactly one quadrilateral. For example, the edge
$p0_j\edash p1_k$ is a side of the quadrilateral  having
$e2_{-j-k-1}$ as the opposite side of $p0_j\edash p1_k$; recall that
$e2_{-j-k-1}$ is in $L_2$ between the two vertices $p2_{-j-k-2}$ and
$p2_{-j-k-1}$, which are, respectively, the vertices of the two
triangles in $S_n$ having $p0_j\edash p1_k$ as a common side. See
Figure \ref{quadpic.fig}.

The quadrilateral is formed by minimizing geodesics joining the
points on $e2_{-j-k-1}$ to the distance-proportional points on
$p0_j\edash p1_k$. In particular, the two edges in $L_0$, $L_1$
joining $p20$, $p12$ and $p0_j$, $p1_k$, respectively, are the other
two sides of the quadrilateral. We denote this quadrilateral by
$q01_{jk}$ and the set of quadrilaterals is
\begin{equation}\label{eqn:reclist}
q01_{jk}, \ \ q12_{jk}, \ \ q20_{jk}, \ \ \ j,k=0,\dots,n-1 \ \
\text{(mod $n$)}\ .
\end{equation}

This concludes our construction of the $2$-cells.  The only new
vertices added are then those in (\ref{eqn:cp2poles}) and the only
new edges are those in the $L_j$'s.

 We now proceed to describe the $3$-cells. It should be pointed
out that, up to now, the cells constructed can be easily verified to
be embedded in $\mathbb C\mathbf P^2$, and there is no intersection
among them in the interior of any cell. As the dimension of the cell
becomes higher, this becomes less clear {\it a priori\/}. We will
show later that the cells do form a tessellation for the $\mathbb
C\mathbf P^2$.

The $3$-cells are divided into two groups. Each of them is in
two-to-one correspondence with the set of edges in $S_n$, or the set
of quadrilaterals. In fact, every quadrilateral is a face of exactly
two $3$-cells in each group.

\begin{figure}[!htbp]
\centerline{ \psfig{figure=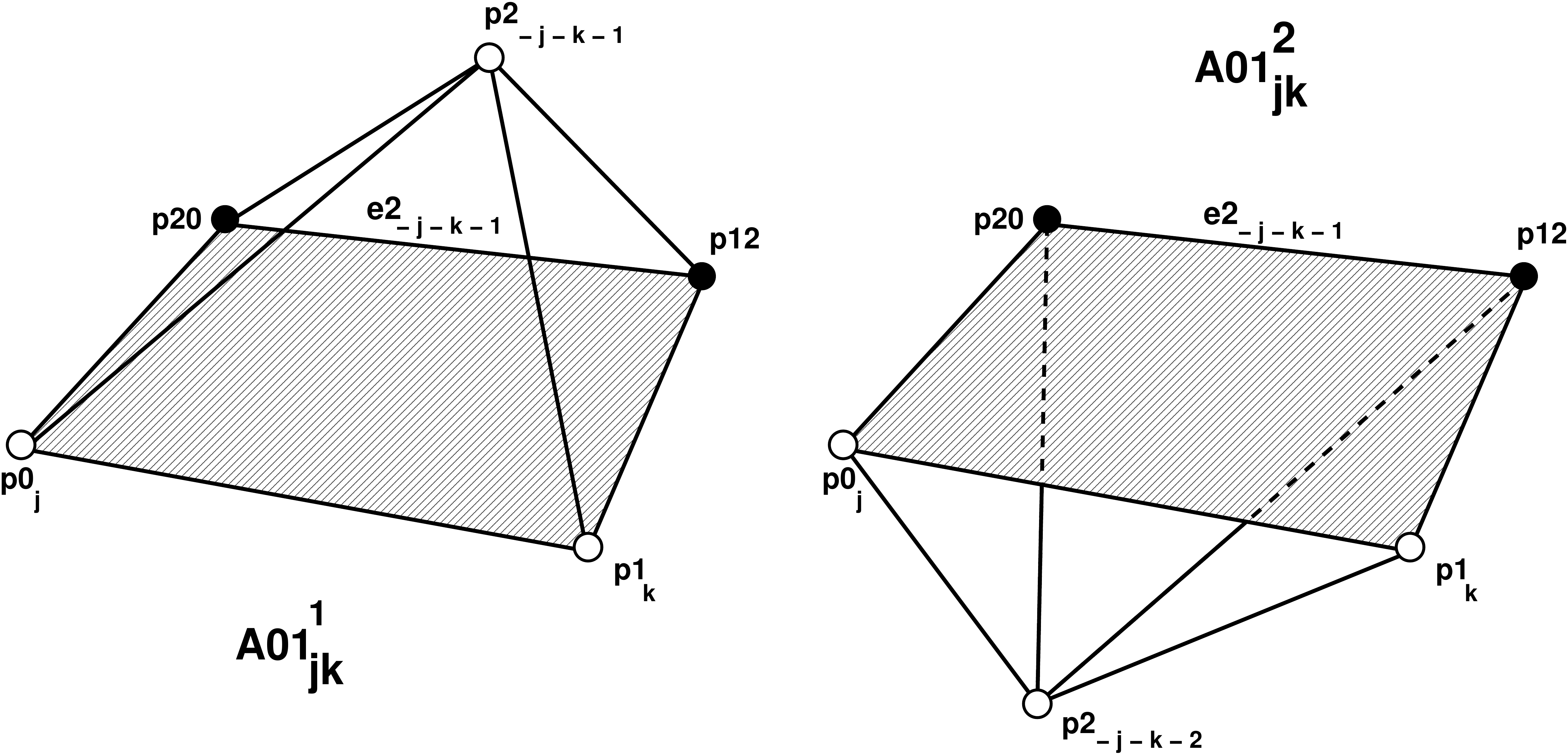,height=6.5cm}}
\caption{}
  \label{groupA3cells.fig}
\end{figure}

In the first group, the two $3$-cells corresponding to, say, the
edge $p0_j\edash p1_k$ are formed by interpolating between
distance-proportional points on $b_{j,k,-j-k-1}$, $b_{j,k,-j-k-2}$
and $f2_{-j-k-1}^{-j-k-1}$, $f2_{-j-k-1}^{-j-k-2}$, respectively, by
minimizing geodesics (see Figure \ref{groupA3cells.fig}).

Clearly, the $3$-cell is a pyramid. Besides the quadrilateral face
$q01_{jk}$, the other four faces are the triangles
$$\{ \ b_{j,k,-j-k-1}, \ \ h12_{k,-j-k-1}, \ \ h20_{-j-k-1,j}, \ \
f2_{-j-k-1}^{-j-k-1} \ \},$$ or
$$\{ \ b_{j,k,-j-k-2}, \ \ h12_{k,-j-k-2}, \ \ h20_{-j-k-2,j}, \ \ f2_{-j-k-2}^{-j-k-1} \ \},$$
respectively. Denote these pyramids by $A01_{jk}^1$, $A01_{jk}^2$,
respectively. We can now list all the $6n^2$ $3$-cells in the first
group:
\begin{eqnarray}\label{eqn:pyrAlist}
&& A01_{jk}^1, \ \ A01_{jk}^2, \nonumber \\
&& A12_{jk}^1, \ \ A12_{jk}^2, \ \ \ j,k=0,\dots,n-1 \
\ \text{(mod $n$)}. \\
&& A20_{jk}^1, \ \ A20_{jk}^2. \nonumber
\end{eqnarray}

\begin{figure}[!htbp]
 \centerline{
\psfig{figure=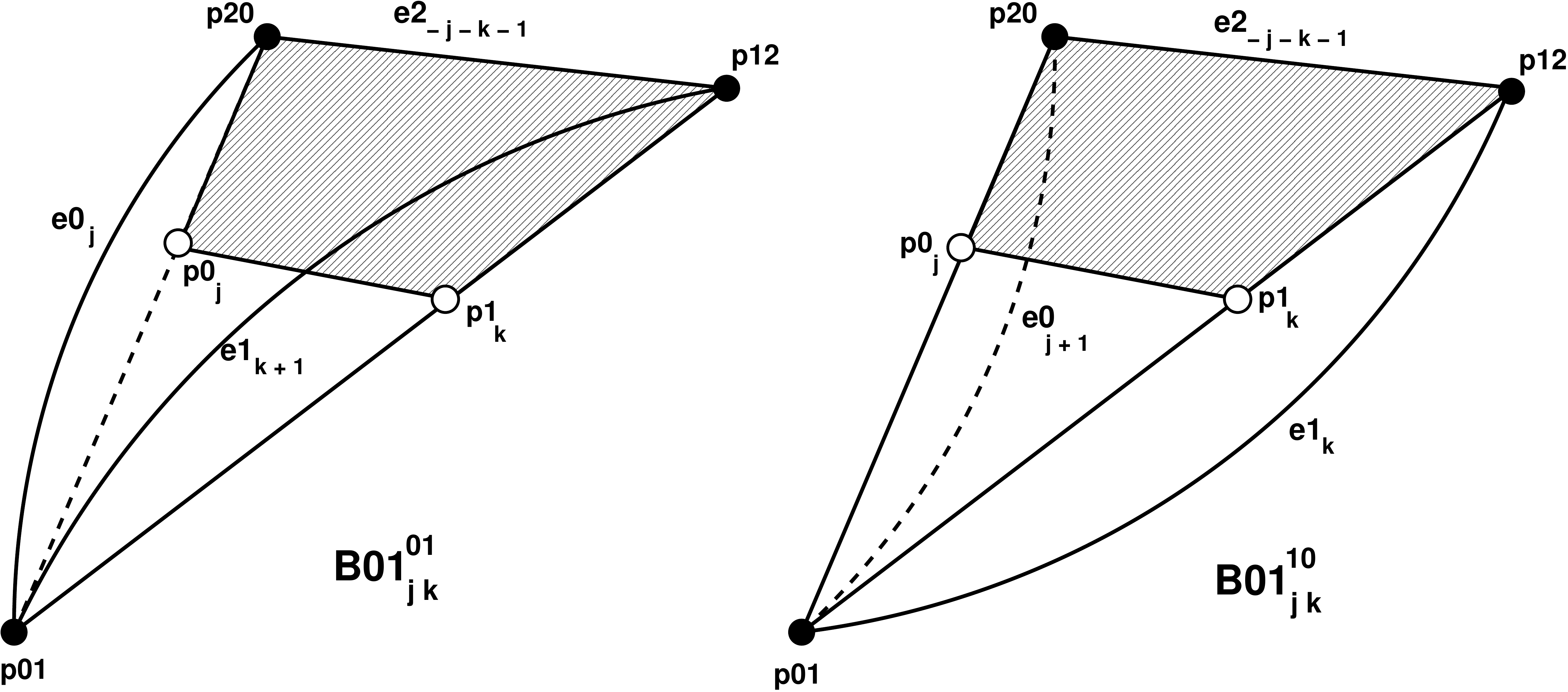,height=6.5cm} }
\caption{}
  \label{groupB3cells.fig}
\end{figure}

In the second group of $3$-cells, the two corresponding to, say
again, $p0_j\edash p1_k$ are formed by minimizing geodesic
interpolation between $h01_{jk}$ and  $g_{j,k+1,-j-k-1}$,
$g_{j+1,k,-j-k-1}$, respectively (see Figure
\ref{groupB3cells.fig}).

Clearly, each $3$-cell is also a pyramid. Besides the quadrilateral
face $q01_{jk}$, the other four faces are the triangles
$$\{ \ g_{j,k+1,-j-k-1}, \ \ f0_j^j, \ \ f1_{k+1}^k, \ \ h01_{jk} \ \}\ ,$$
or
$$\{ \ g_{j+1,k,-j-k-1}, \ \ f0_{j+1}^j, \ \ f1_k^k, \ \ h01_{jk} \ \}\ ,$$
respectively. Notice that, unlike the first group, whose two
pyramids share only the quadrilateral face, these two pyramids share
both the quadrilateral face $q01_{jk}$ and the triangular face
$h01_{jk}$. Denote these pyramids by $B01_{jk}^{01}$,
$B01_{jk}^{10}$, respectively. The list of all the $6n^2$ $3$-cells
in the second group then is:
\begin{eqnarray}\label{eqn:pyrBlist}
&& B01_{jk}^{01}, \ \ B01_{jk}^{10}, \nonumber \\
&& B12_{jk}^{01}, \ \ B12_{jk}^{10}, \ \ \ j,k=0,\dots,n-1 \
\ \text{(mod $n$)}. \\
&& B20_{jk}^{01}, \ \ B20_{jk}^{10}\ . \nonumber
\end{eqnarray}

 We are now ready to tessellate the $\mathbb C\mathbf P^2$ by
$4$-cells. Each $4$-cell is bounded by four pyramids, two from each
of the groups (\ref{eqn:pyrAlist}) and (\ref{eqn:pyrBlist}); in
fact, two from one determine the two from the other. Since every
$3$-cell should be the face of exactly two $4$-cells, it follows
that there are in all $6n^2$ $4$-cells. We illustrate one of them as
follows.

Start with $A01_{jk}^1$ in (\ref{eqn:pyrAlist}). The other pyramid
from (\ref{eqn:pyrAlist}) is either $A12_{k,-j-k-1}^1$ or
$A20_{-j-k-1,j}^1$, as these are the only other two pyramids in
(\ref{eqn:pyrAlist}) sharing the triangular face $b_{j,k,-j-k-1}$
with $A01_{jk}^1$. If, say, we pick $A12_{k,-j-k-1}^1$, then it is
easy to see that the two pyramids from (\ref{eqn:pyrBlist}) must be
$B01_{jk}^{01}$ and $B12_{k,-j-k-1}^{10}$, in order to have the
quadrilateral faces $q01_{jk}$ and $q12_{k,-j-k-1}$ shared, and for
the two to have the triangular face from (\ref{eqn:trilist4}) in
common. Therefore this $4$-cell is bounded by the following four
pyramids:
\begin{equation}\label{eqn:3sphere}
\{ \ A01_{jk}^1, \ \ A12_{k,-j-k-1}^1, \ \ B01_{jk}^{01}, \ \
B12_{k,-j-k-1}^{10} \ \}
\end{equation}

\begin{figure}[!htbp]
\centerline{ \psfig{figure=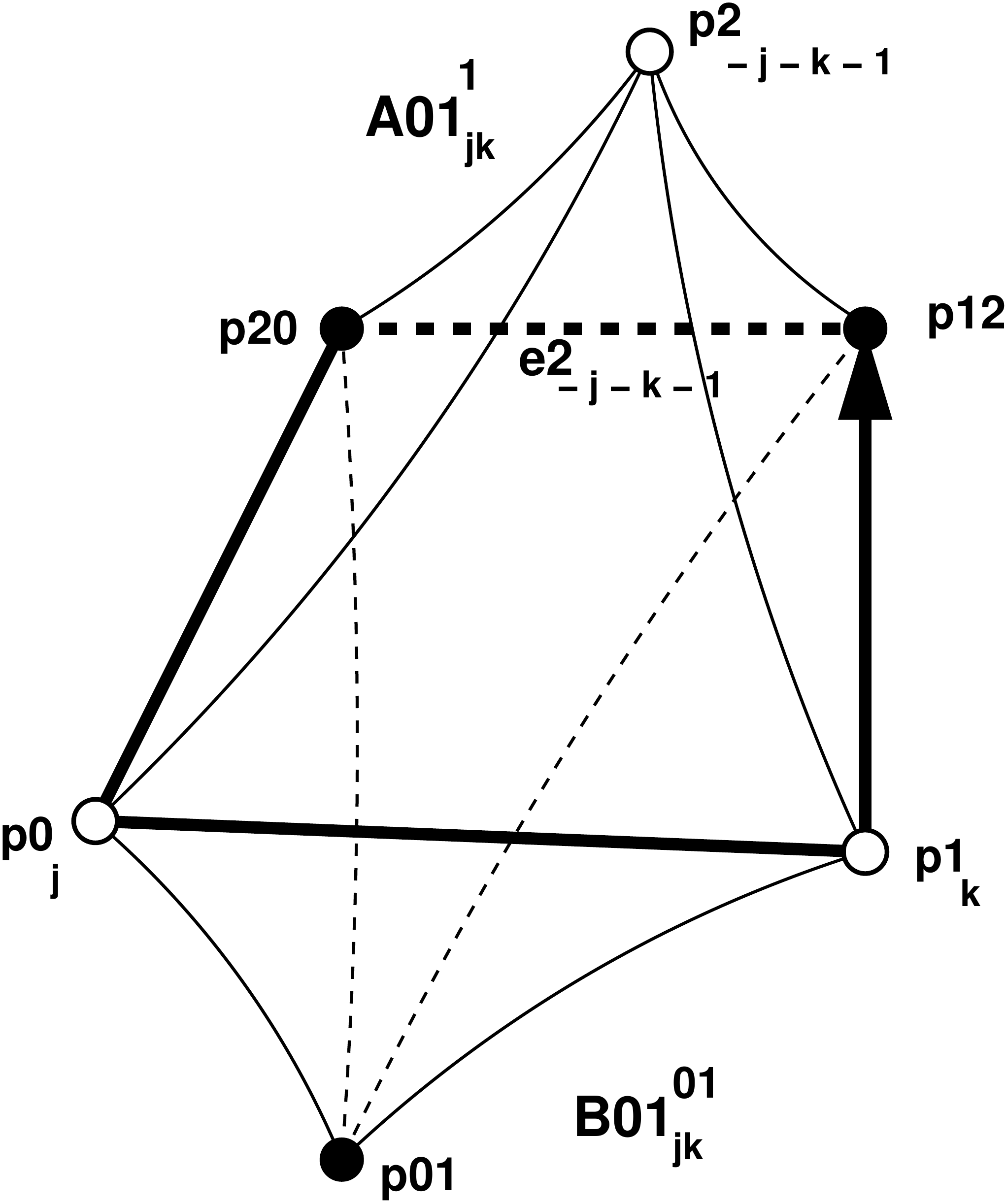,height=5.5cm}
\hspace{0.5in}    \psfig{figure=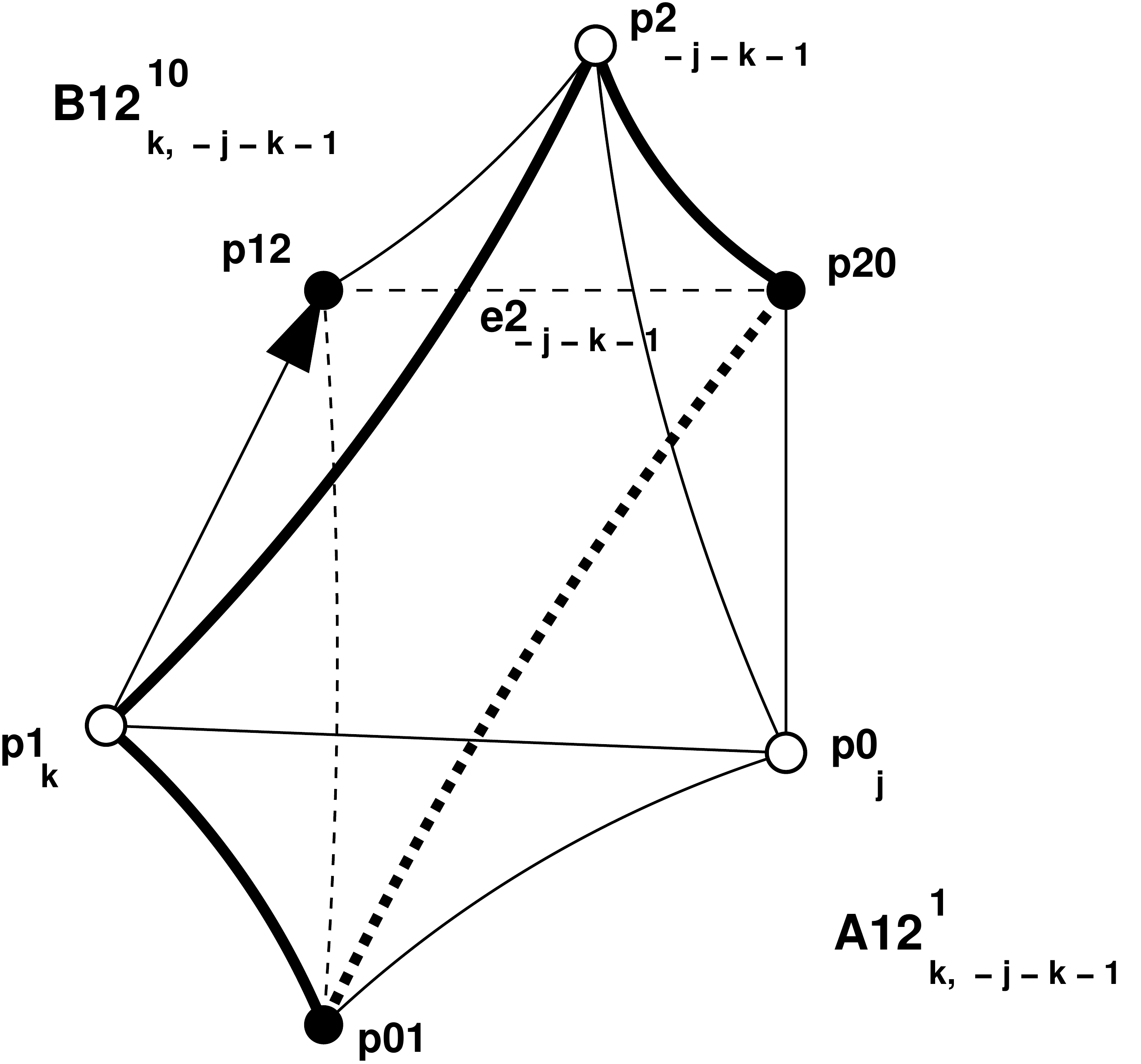,height=5.5cm} }
\caption{}
  \label{S3tess.fig}
\end{figure}

As illustrated in Figure \ref{S3tess.fig}, the pyramids in
(\ref{eqn:3sphere}) indeed form a tessellation for a 3-sphere, at
least combinatorially.

From the above, it is easy now to list all the $4$-cells in terms of
their boundary pyramids:
\begin{eqnarray}\label{eqn:4celist}
&&\{A01_{jk}^1, \ A12_{k,-j-k-1}^1, \ B01_{jk}^{01},
\ B12_{k,-j-k-1}^{10}\}, \nonumber \\
&&\{A01_{jk}^2, \ A12_{k,-j-k-2}^2, \ B01_{jk}^{10},
\ B12_{k,-j-k-2}^{01}\}, \nonumber \\
&& \{A12_{jk}^1, \ A20_{k,-j-k-1}^1, \ B12_{jk}^{01},
\ B20_{k,-j-k-1}^{10}\}, \\
&& \{A12_{jk}^2, \ A20_{k,-j-k-2}^2, \ B12_{jk}^{10},
\ B20_{k,-j-k-2}^{01}\}, \nonumber \\
&& \{A20_{jk}^1, \ A01_{k,-j-k-1}^1, \ B20_{jk}^{01}, \
B01_{k,-j-k-1}^{10}\}, \nonumber \\
&& \{A20_{jk}^2, \ A01_{k,-j-k-2}^2, \ B20_{jk}^{10}, \
B01_{k,-j-k-2}^{01}\}. \nonumber
\end{eqnarray}
for $j,k=0,\dots,n-1 \ \ \text{(mod $n$)}$. Combinatorially, these
$6n^2$ $4$-cells together form a simplicial $4$-manifold. Combining
this with the numbers of vertices, edges, $2$-cells, and $3$-cells
we have obtained before, we find its Euler characteristic number to
be
$$(3n+3)-(3n^2+9n)+(2n^2+6n+3n^2+n^2+3n^2)-12n^2+6n^2=3 \ ,$$
which is the Euler characteristic of $\mathbb C\mathbf P^2$.
However, as we pointed out earlier, to show this is really a
tessellation of the $\mathbb C\mathbf P^2$, one needs to verify that
all the $4$-cells are embedded and that there is no intersection
among them at any of their interior points. We now confirm this.

For any fixed point $p=[0,z_1,z_2]\in L_0$ let $L_{0,p}$ be the
projective line joining $p12$ and $p$. Then $\mathbb C\mathbf
P^2=\bigcup_{p\in L_0}L_{0,p}$; the union is disjoint except that
all the $L_{0,p}$'s intersect at the single point $p12$. It is easy
to verify, (i) if $p\notin S_n$, then $L_{0,p}$ intersects $S_n$ at
exactly $n$ different points in a similar position to those in
(\ref{eqn:cp1roots}) on $\mathbb C\mathbf P^1$, and (ii) if $p\in
S_n$ then $p$ is the only intersection of $L_{0,p}$ and $S_n$.

\begin{figure}[!htbp]
\centerline{    \psfig{figure=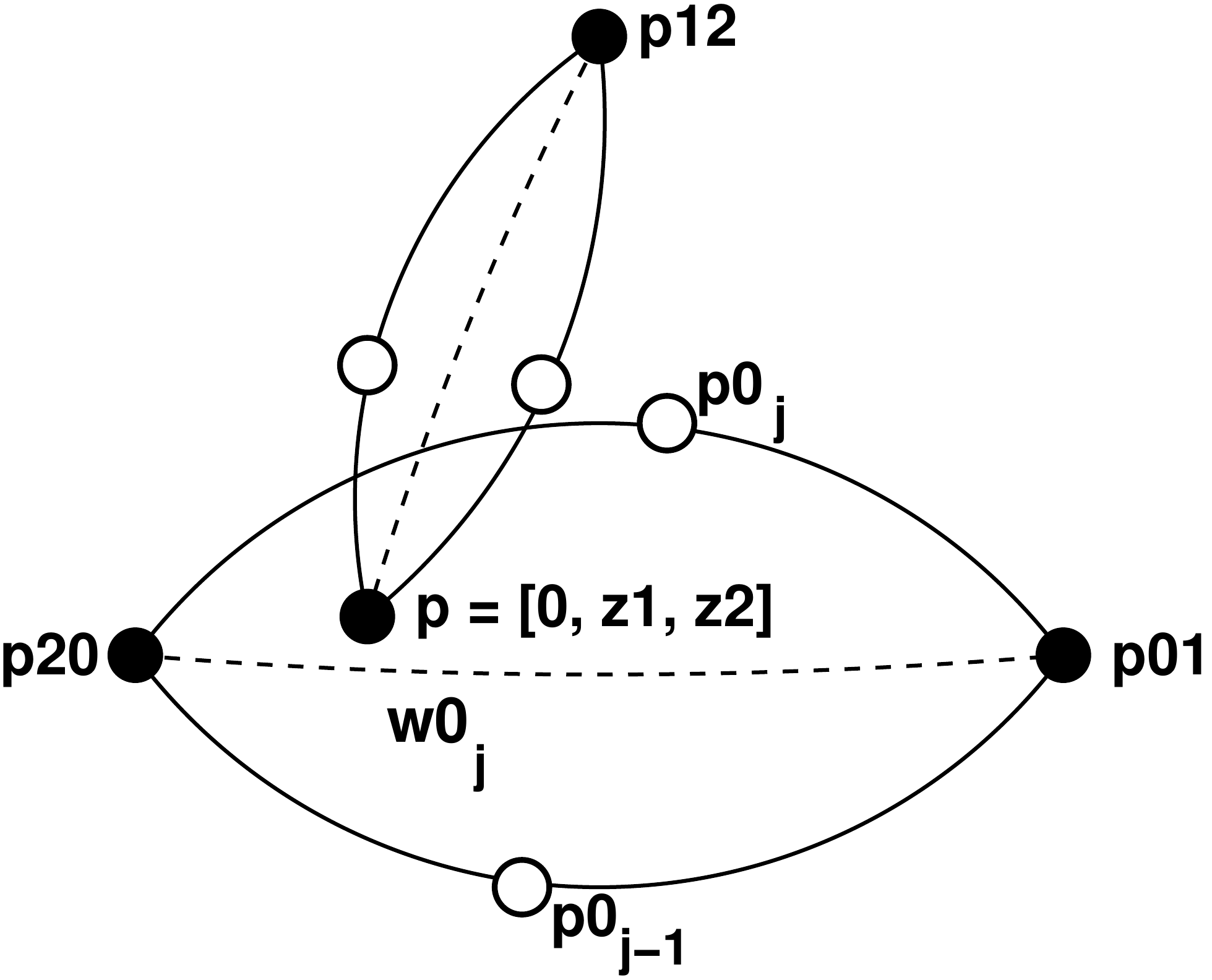,height=5.5cm}
\hspace{.1in}  \psfig{figure=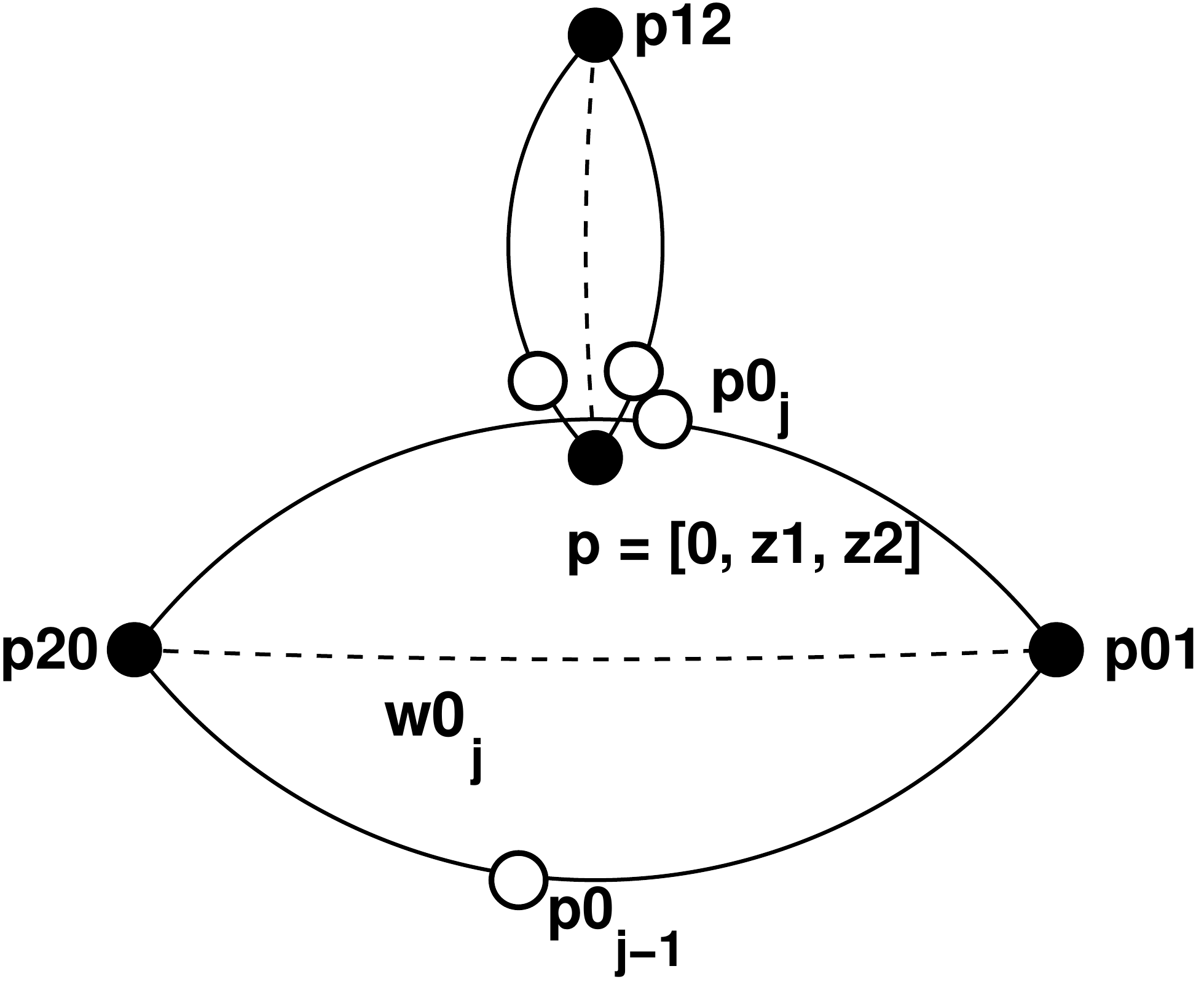,height=5.5cm} }
  \caption[]{}
  \label{sweep1.fig}
\end{figure}

For $p\notin S_n$, we triangulate $L_{0,p}$ similarly to $\mathbb
C\mathbf P^1$, using the points $p12$, $p$ (corresponding to the
vertices in (\ref{eqn:cp1poles})), and the $n$ intersections with
$S_n$ (see Figure \ref{sweep1.fig}).

 Notice that for $p\in S_n$,
although we do not have the triangulation, there are $n$ well
defined paths from $p12$ to $p$ that are obtained as limits of the
paths on $L_{0,q}$, for $q$ near $p\in S_n$ on $L_0$, corresponding to the
edges joining $p0$ in (\ref{eqn:cp1poles}) and the $p_k$'s in
(\ref{eqn:cp1roots}).  Also notice, in particular, that
$L_{0,p01}=L_1$ and $L_{0,p20}=L_2$.

Let
$$wl_k=fl\,{}_k^{k-1}\cup fl\,{}_k^k \ , \ \ \ l=0,1,2; \ \ k=0,\dots,n-1
\ \text{(mod $n$)} \ .$$

For given $j$ and $k$, as $p$ varies on $w0_j$, it is easy to see
that we get a continuous family of regions $w1_{k,p}$ in $L_{0,p}$
with $w1_{k,p01}=w1_k$. The union of this family of regions then
clearly forms an embedded $4$-cell in the $\mathbb C\mathbf P^2$,
which we will denote by $W_{j,k,-j-k}$. It is also clear that
$\mathbb C\mathbf P^2=\bigcup_{j,k}W_{j,k,-j-k}$ and there is no
intersection between different $W_{j,k,-j-k}$'s at any of their
interior points.

\begin{figure}[!htbp]
\centerline{\psfig{figure=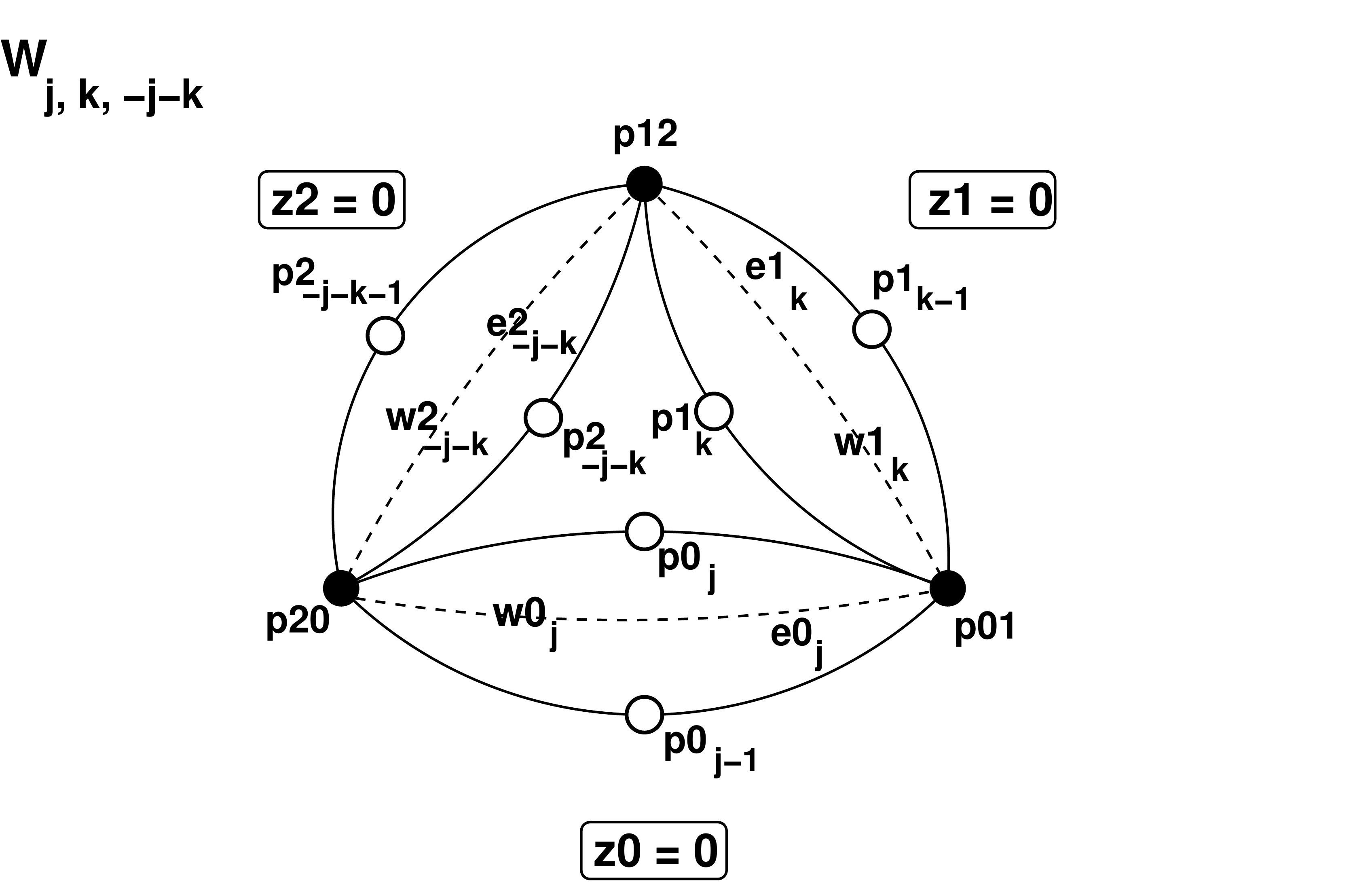,height=8cm} } % 5.5cm
  \caption[]{}
  \label{sweep2.fig}
\end{figure}

For clarity and later convenience, we write down the following
explicit parametrization for $W_{j,k,-j-k}$:
\begin{equation}\label{eqn:sweepcell}
[\cos s\,e^{i\beta},\cos r\sin s,\sin r\sin s\,e^{i\alpha}]
\end{equation}
with \ $0\le r,s\le\pi/2$, \ $(2j-1)\pi/n\le\alpha\le (2j+1)\pi/n$ \
and, if we denote  \ $\arg(-\cos^nr-\sin^nr\,e^{in\alpha})$ \ by
$a(r,\alpha)$ with $0\le a(r,\alpha)\le 2\pi$, then
\begin{equation*}
\frac{a(r,\alpha)+2(j+k-1)\pi}{n}\le\beta\le
\frac{a(r,\alpha)+2(j+k)\pi}{n} \ .
\end{equation*}
We see in particular that $w1_{k,p20}=w2_{-j-k}$.  See Figure
\ref{sweep2.fig}.

From (\ref{eqn:sweepcell}), it follows that, in a way similar to the
above, $W_{j,k,-j-k}$ can also be described as a union of
$w2_{-j-k,p}$ over $p\in w1_k$, or a union of $w0_{j,p}$ over $p\in
w2_{-j-k}$. Therefore the boundary of $W_{j,k,-j-k}$ is tessellated
by twelve $3$-cells; each of them is the union of one of the two
lower half boundary edges of a $w$-region over one of the two
triangles in the corresponding base region. It is easy to see that
these $3$-cells are in fact exactly the following twelve pyramids in
(\ref{eqn:pyrAlist}):
\begin{equation}\label{eqn:sweepboundary}
\begin{array}{llll}
 A12_{k-1,-j-k}^1\ , & A12_{k-1,-j-k}^2 \ ,& A12_{k,-j-k-1}^1 \ ,
& A12_{k,-j-k-1}^2 \ ,\\[.05in]
A20_{-j-k-1,j}^1 \ , & A20_{-j-k-1,j}^2 \ ,& A20_{-j-k,j-1}^1 \ ,& A20_{-j-k,j-1}^2  \ ,\\[.05in]
A01_{j-1,k}^1 \ ,& A01_{j-1,k}^2 \ ,& A01_{j,k-1}^1 \ ,&
A01_{j,k-1}^2 \ .
\end{array}
\end{equation}

It is also easy to see there are six $3$-cells contained inside
$W_{j,k,-j-k}$; each of them is the union of one of the two
triangles in a $w$-region over the middle edge of the corresponding
base region. These $3$-cells are the following six pyramids in
(\ref{eqn:pyrBlist}):
\begin{equation}\label{eqn:sweepinside}
\begin{array}{ll}
 B12_{k-1,-j-k}^{10} \ , & B12_{k,-j-k-1}^{01}  \\[0.05in]
 B20_{-j-k-1,j}^{10} \ , & B20_{-j-k,j-1}^{01} \ , \\[0.05in]
 B01_{j-1,k}^{10} \ , & B01_{j,k-1}^{01} \ .
 \end{array}
\end{equation}
These $3$-cells divide $W_{j,k,-j-k}$ into six $4$-cells; they are
the following six in (\ref{eqn:4celist}):
\begin{equation}\label{eqn:Wdivide}
\begin{array}{llll}
\{\ A01_{j,k-1}^1\ , & A12_{k-1,-j-k}^1\ , & B01_{j,k-1}^{01}\ ,
 & B12_{k-1,-j-k}^{10}\ \},  \\[0.05in]
\{\ A12_{k-1,-j-k}^2\ , & A20_{-j-k,j-1}^2\ , & B12_{k-1,-j-k}^{10}\
,
 & B20_{-j-k,j-1}^{01}\ \}\ ,  \\[0.05in]
\{\ A20_{-j-k,j-1}^1\ , & A01_{j-1,k}^1\ , & B20_{-j-k,j-1}^{01}\ ,
 & B01_{j-1,k}^{10}\ \}\ , \\[0.05in]
\{\ A01_{j-1,k}^2\ ,& A12_{k,-j-k-1}^2\ , & B01_{j-1,k}^{10}\ ,
 &  B12_{k,-j-k-1}^{01}\ \}\ ,  \\[0.05in]
\{\ A12_{k,-j-k-1}^1\ , & A20_{-j-k-1,j}^1\ , & B12_{k,-j-k-1}^{01}\
,
 &  B20_{-j-k-1,j}^{10}\ \}\ , \\[0.05in]
\{\ A20_{-j-k-1,j}^2\ , & A01_{j,k-1}^2\ , & B20_{-j-k-1,j}^{10}\ ,
 &B01_{j,k-1}^{01}\ \}\ .
\end{array}
\end{equation}

\begin{figure}[!htbp]
 \centerline{
\psfig{figure=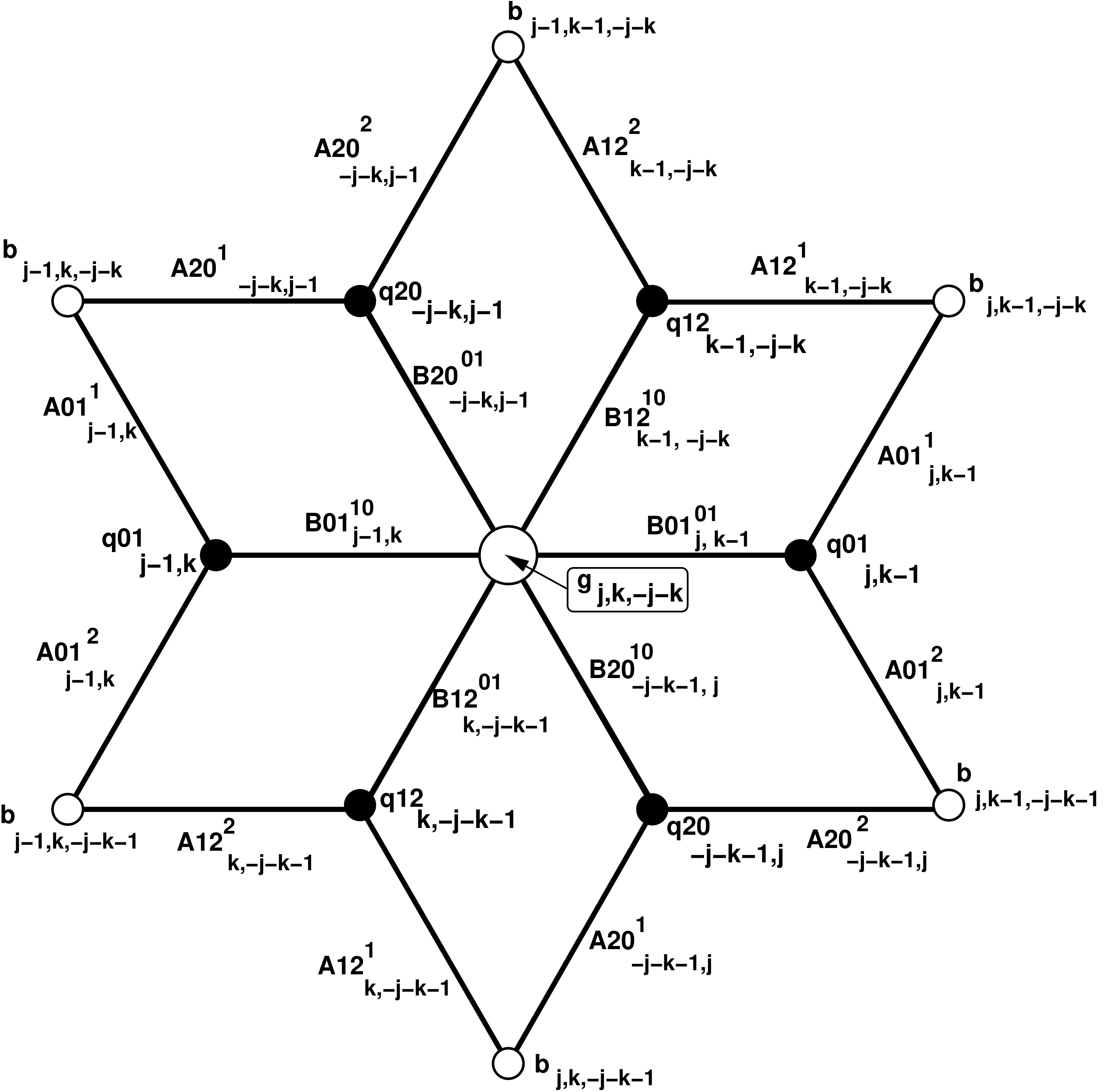,height=10.0cm} }
  \caption[]{}
   \label{hexagon.fig}
\end{figure}

The structure of (\ref{eqn:sweepboundary}), (\ref{eqn:sweepinside})
and (\ref{eqn:Wdivide}) together can be illustrated by the diagram
in Figure  \ref{hexagon.fig}.

Now it is clear that the $\mathbb C\mathbf P^2$ is well tessellated.

{\bf Remark.} If all we need is a tessellation of the $\mathbb
C\mathbf P^2$ (and hence $F_n$), then the triangles in
(\ref{eqn:trilist4}) and the $3$-cells in (\ref{eqn:pyrBlist}) are
not needed. The pyramids in (\ref{eqn:pyrAlist}) are paired into
$3n^2$ octahedrons, and the $4$-cells of the tessellation are
precisely the $n^2\,$ $W_{j,k,-j-k}$'s.   However, when this
tessellation is lifted to $F_n$, it is not $\Gamma_n$-invariant.

\bigskip

\section{Tessellation of $F_n$}

Through the $n$-fold regular branched covering
(\ref{eqn:branchcover}), the tessellation of $\mathbb C\mathbf P^2$
in $\S$\ref{section:tesscp2} now lifts to a well defined
tessellation for the $F_n$. The numbers of vertices, edges,
$2$-cells, $3$-cells and $4$-cells are as indicated in the
introduction. In this section, we examine this tessellation more
closely and show that it is $\Gamma_n$-invariant. Recall that
$\Gamma_n$ is the group of the isomorphisms of $F_n$ induced from
permuting and/or phase multiplying the homogeneous coordinates of
$\mathbb C\mathbf P^3$.

The intersection of $F_n$ with each of the projective planes $z_k=0,
\ \ k=0,1,2,3$, is the $S_n$ in that plane, triangulated as
described in $\S$\ref{section:tesscurve}. The tessellation of the
$F_n$ is an extension of the triangulations on these four $S_n$'s.
In fact, the four $S_n$'s contain all the vertices, edges, and
triangles lifted from those in (\ref{eqn:trilist2}) and
(\ref{eqn:trilist1}). The other $4n^3$ triangles are lifted from
(\ref{eqn:trilist3}) and (\ref{eqn:trilist4}) and are  characterized
by the fact that for each of them, the three edges lie on three
distinct $S_n$'s. Notice then that for any three $S_n$'s of the
four, any three different pairwise intersections are vertices of a
unique triangle lifted from (\ref{eqn:trilist3}) or
(\ref{eqn:trilist4}).

The formation of the quadrilaterals can be described as follows.
Start with any edge on one of the four $S_n$'s, say, the one in the
projective plane $z_0=0$; its two end points, denoted by $q_1$ and
$q_2$, must then also lie in two other distinct projective planes,
say, $z_1=0$ and $z_2=0$, respectively. Then there are $n$ distinct
edges on $z_2=0$ joining $q_2$ and the $n$ distinct intersections of
the projective planes $z_2=0$ and $z_3=0$. Any one of these edges
plus $\overline{q_1q_2}$, the edge we started with, form two
adjacent sides of a unique quadrilateral. Hence one sees that there
are in all $3n^3$ quadrilaterals.

\begin{figure}[!htbp]
 \centerline{
\psfig{figure=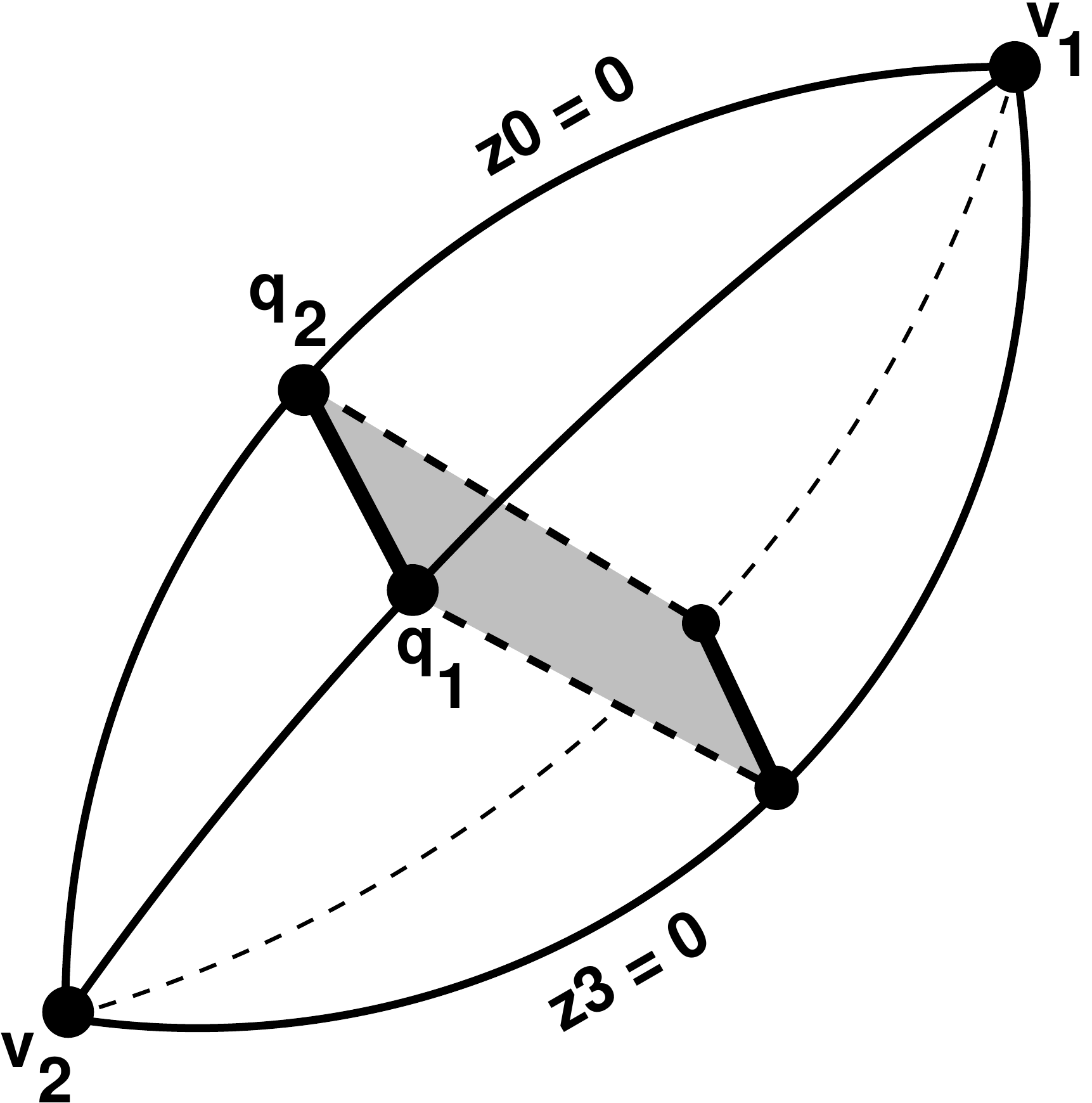,height=6.0cm} }
  \caption[]{}
   \label{v1v2.fig}
\end{figure}

For two opposite sides, say, lying in the projective planes $z_0=0$
and $z_3=0$, respectively, of a given quadrilateral,  as in the
example above, there are exactly two vertices, $v_1,v_2$, in the
intersection of $z_0=0$ and $z_3=0$ that are the opposite vertices
of the given edges in triangles lying in $z_0=0$ and $z_3=0$,
respectively (see Figure \ref{v1v2.fig}). Each of these two vertices
forms a pyramid with the quadrilateral. Notice that if one starts
with the other pair of opposite sides of the quadrilateral, the two
vertices will be different. One sees that there are in all $12n^3$
pyramids.

Finally, every $4$-cell is bounded by four pyramids, and each
pyramid is shared by two $4$-cells, thus there are $6n^3$ $4$-cells.

From the description above, one can see that if, instead of
(\ref{eqn:branchcover}) which is induced from the projection
$[z_0,z_1,z_2,z_3]\mapsto [z_0,z_1,z_2]$, we use the branched
covering induced from, say, $[z_0,z_1,z_2,z_3]\mapsto
[z_1,z_2,z_3]$, the lifted tessellation will be the same. This,
combined with the invariance for $S_n$ demonstrated in
$\S$\ref{section:tesscurve}, shows that the tessellation of the
$F_n$ is $\Gamma_n$-invariant and the action of $\Gamma_n$ is
transitive on the set of $4$-cells.

The list of all the vertices in the tessellation is:
$$p01_{k}:=[0,0,1,e^{i(\pi+2k\pi)/n}], \ p02_{k}:=[0,1,0,e^{i(\pi+2k\pi)/n}],$$
$$p03_{k}:=[0,1,e^{i(\pi+2k\pi)/n},0], \
p12_{k}:=[1,0,0,e^{i(\pi+2k\pi)/n}],$$
$$p13_{k}:=[1,0,e^{i(\pi+2k\pi)/n},0], \
p23_{k}:=[1,e^{i(\pi+2k\pi)/n},0,0],$$ for $k=0,\dots,n-1$.

\begin{figure}[!htbp]
\centerline{ \psfig{figure=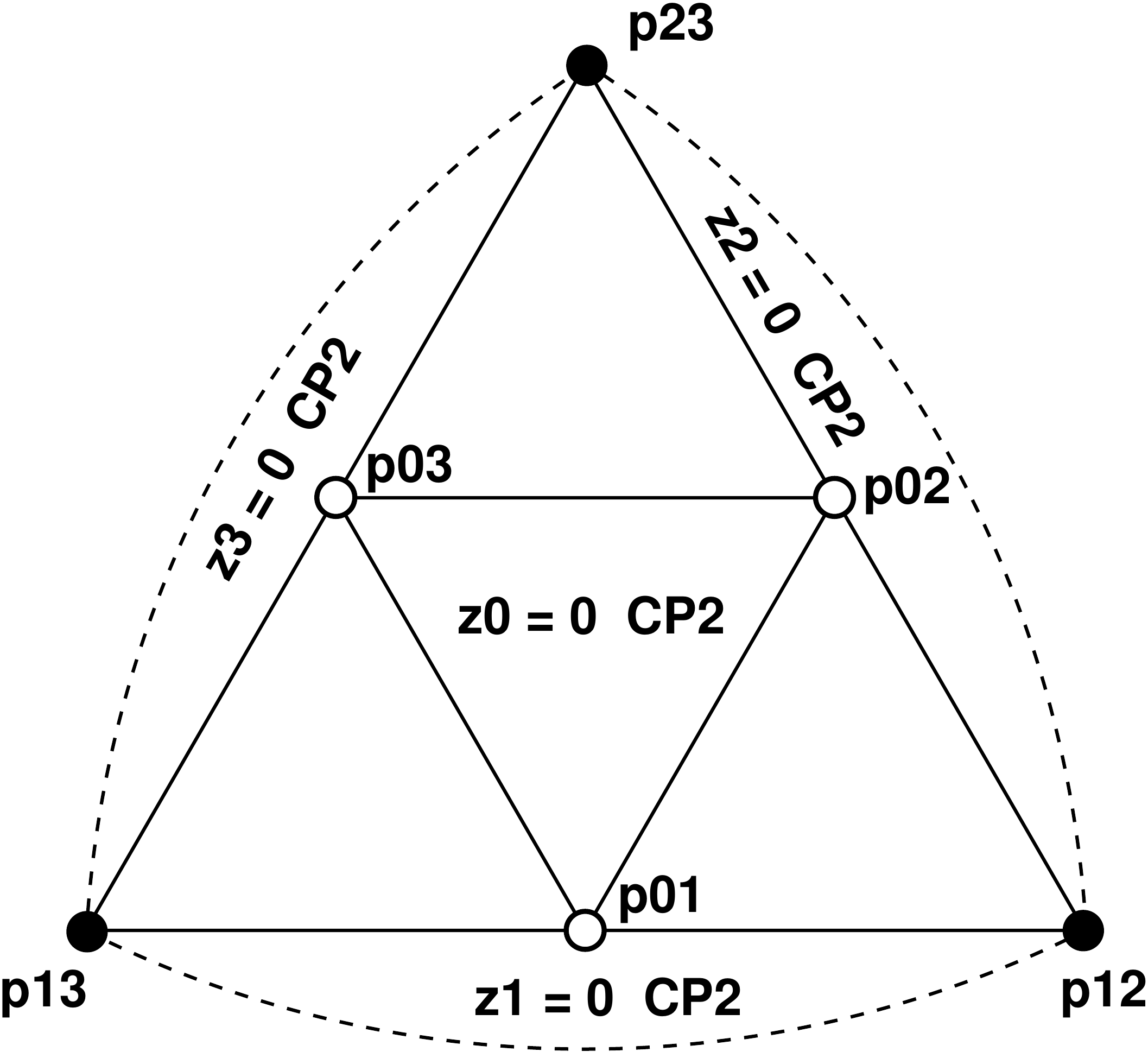,height=5.5cm}
\hspace{0.5in}    \psfig{figure=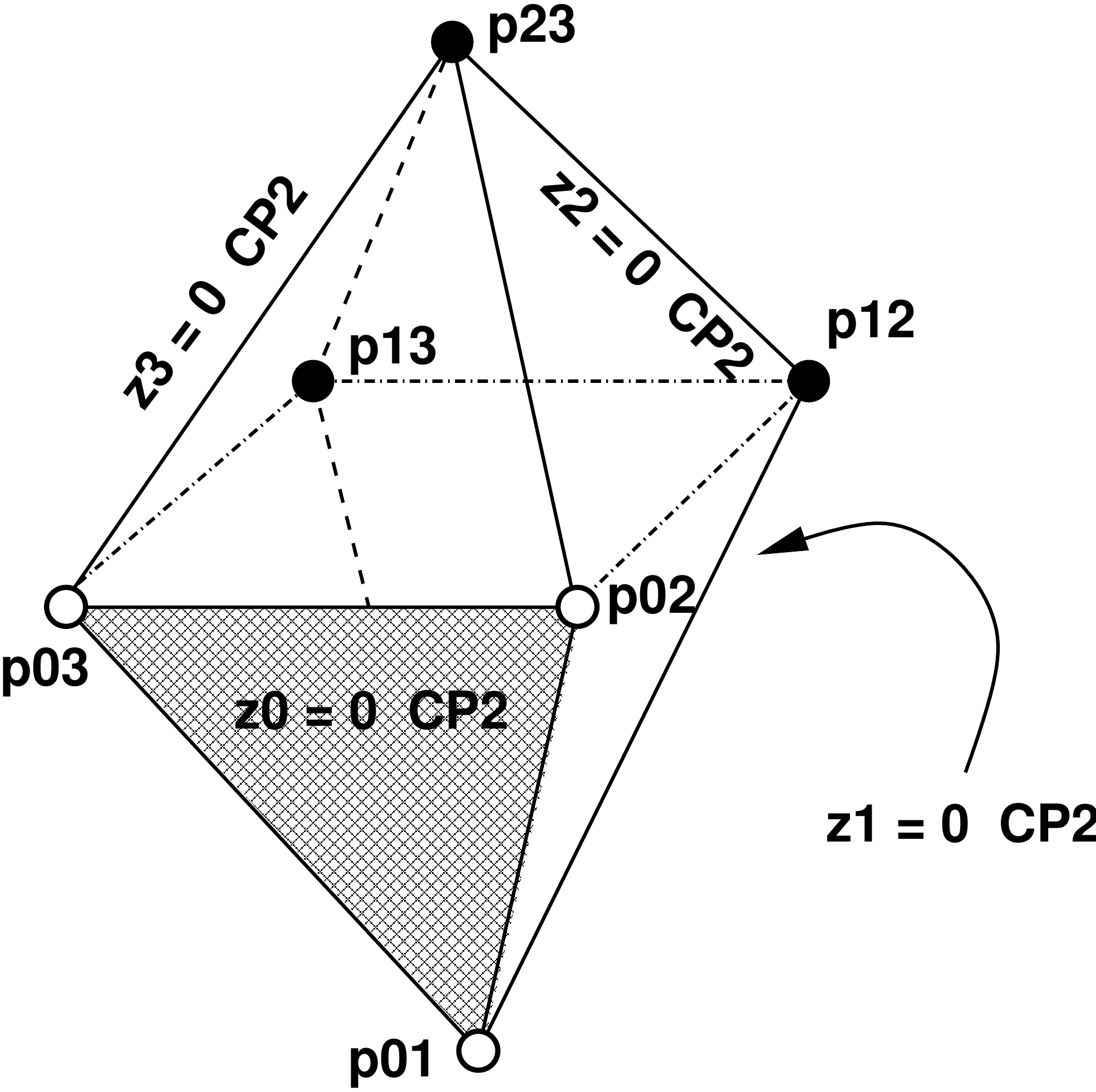,height=5.5cm} }
\caption{}
  \label{CP3-pts.fig}
\end{figure}

One can then list the edges, 2-cells, 3-cells, and 4-cells in terms
of the vertices. We now write down a few lists of edges and 2-cells for
illustration (see Figure \ref{CP3-pts.fig}).

Every edge lies on one of the four $\mathbb C\mathbf P^2$'s defined
by $z_k=0$. For example, the $3n^2$ edges lying on $z_0=0$ are:
$$\{p01_i,\, p02_j\}, \ \ \{p01_i,\, p03_j\}, \ \ \{p02_i,\,
p03_j\},$$ $i,j=0,\dots,n-1$.

The 2-cells are divided into three groups: triangles each of which
lies on one of the four $\mathbb C\mathbf P^2$'s defined by $z_k=0$;
triangles each of which has three sides on three different $\mathbb
C\mathbf P^2$'s; rectangles.

For example, the $2n^2$ triangles lying on $z_0=0$ are:
$$\{p01_i,\, p02_j,\, p03_k\}, \ \ i-j+k=0  \  \ \text{or} \ -1 \ \
(\text{mod} \ n).$$

The $n^3$ triangles whose edges lie on three different $\mathbb
C\mathbf P^2$'s, $z_1=0,\ z_2=0, \ z_3=0,$ respectively, are:
$$\{p12_i,\, p13_j,\, p23_k\}, \ \ \ i,j,k=0,\dots,n-1.$$

The rectangles can be divided into three groups: each of them has
one edge lying on the $\mathbb C\mathbf P^2$ labeled by $z_0=0$ and an
opposite edge on $z_m=0, \ m=1,2,3$. The group with $m=1$, for
example, contains the following $n^3$ rectangles:
$$\{p02_i,\, p03_j,\, p13_k,\, p12_l\}, \ \ \ i-j+k-l=0 \ \ \ (\text{mod} \
n).$$

\begin{figure}[!htbp]
 \centerline{
\psfig{figure=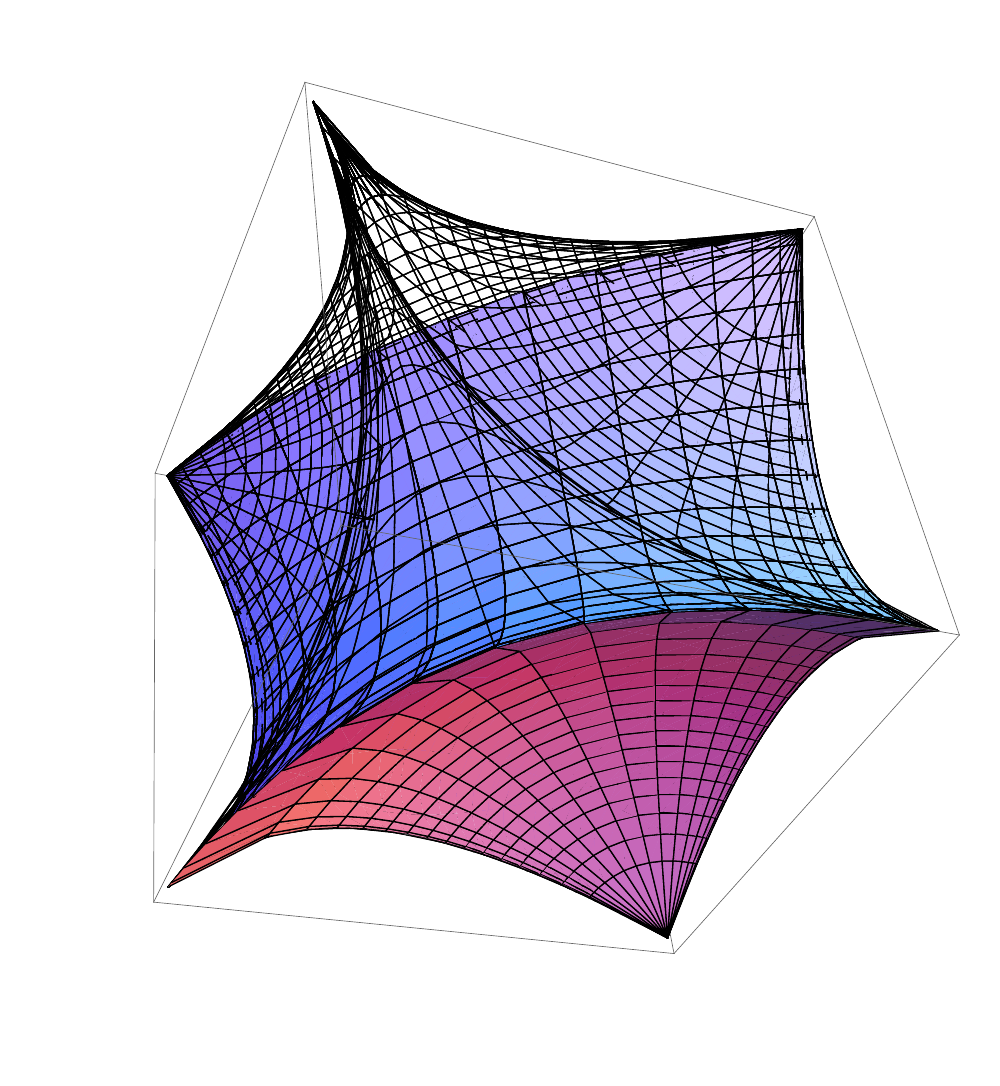,height=8.25cm} \hspace{0.25in}
\psfig{figure=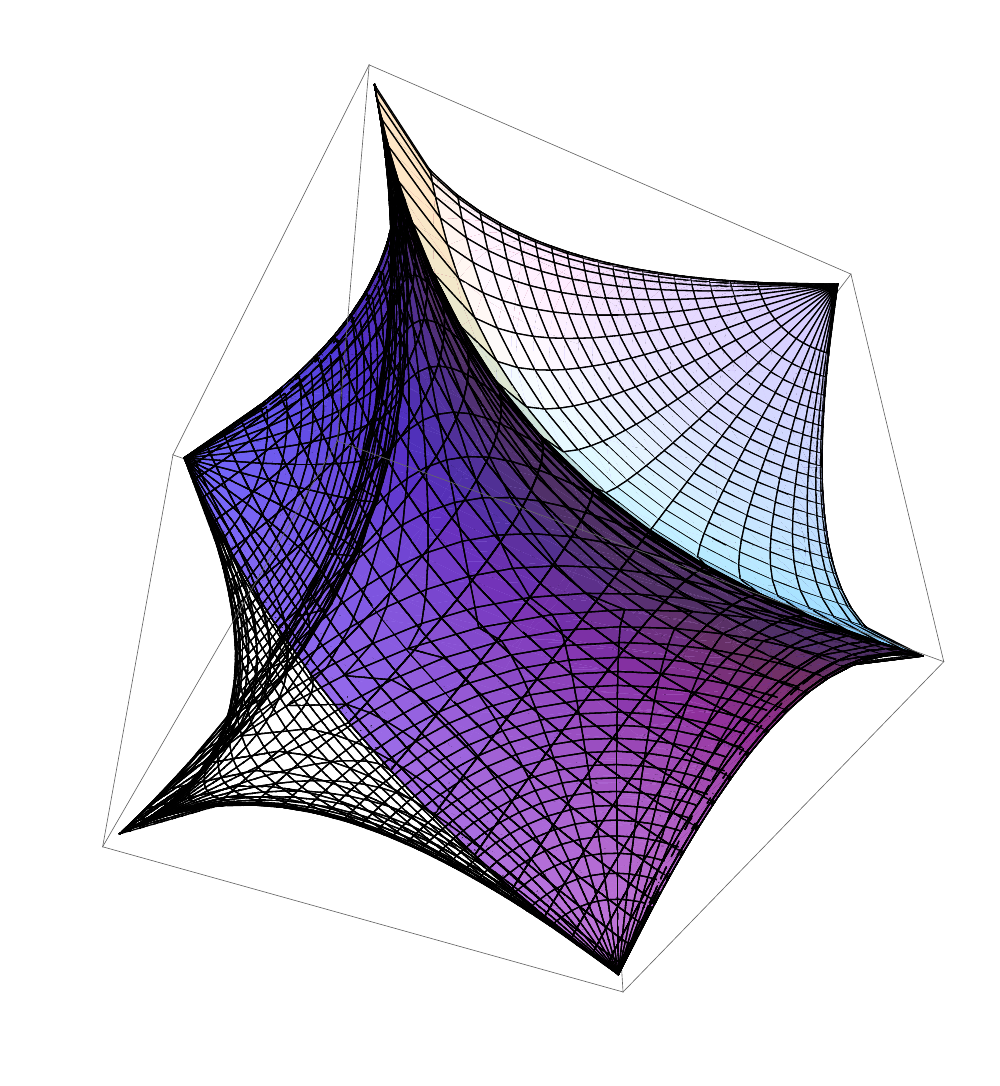,height=8.25cm} } \centerline{
\hspace{1in} (a) \hfill (b) \hspace{1in} }
 \centerline{
\psfig{figure=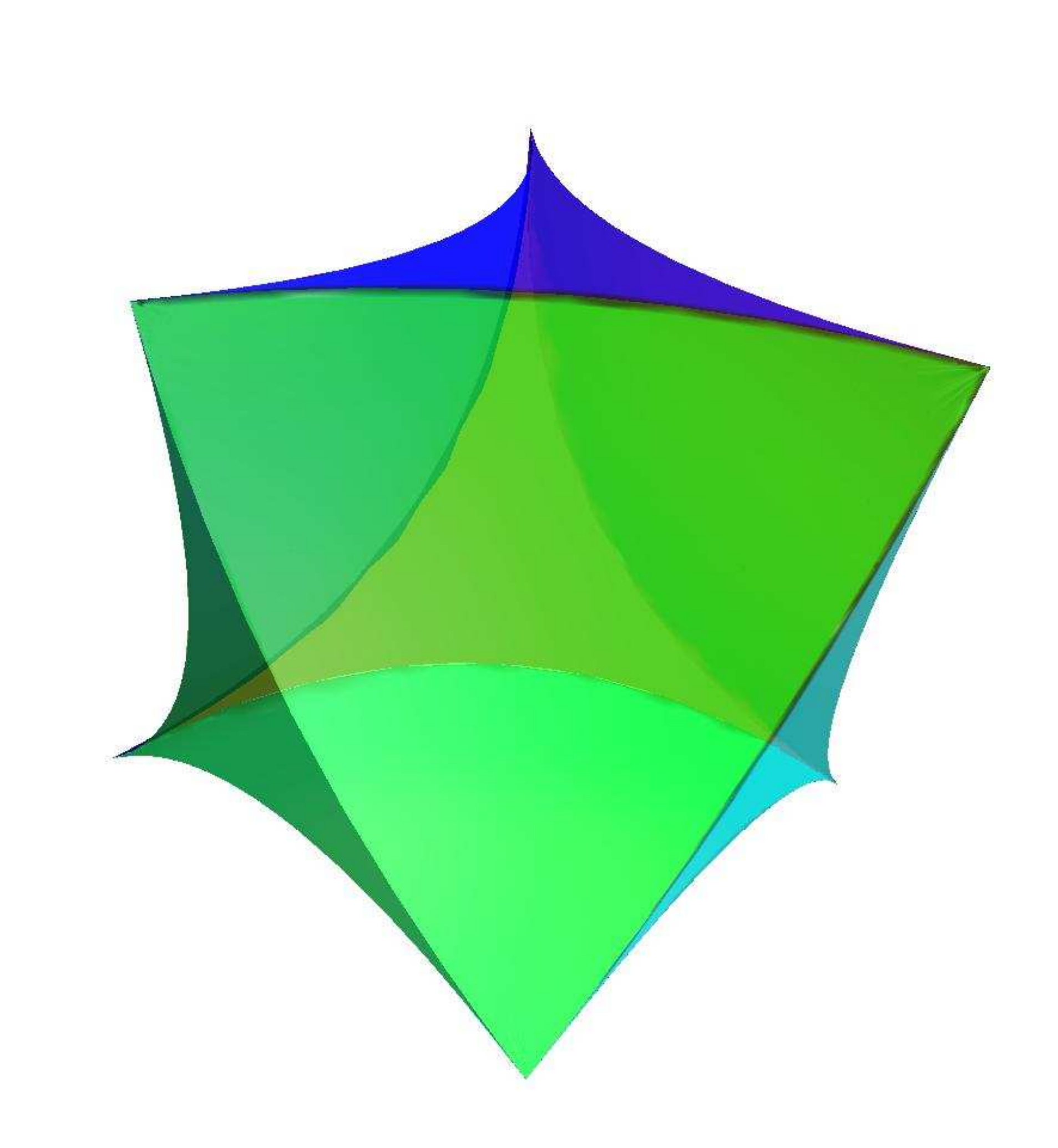,height=10.0cm} }
\centerline{ \hfill
(c) \hfill  }
  \caption[]{}
   \label{4cell.fig}
\end{figure}

Figure \ref{4cell.fig} is an actual image of a generic 4-cell  using
an explicit embedding of $\mathbb C\mathbf P^3$ into
$\mathbb{R}^{16}$ (see, e.g., \cite{ajhjps06}).
 (a) and (b) depict the 3-balls that are the upper
and lower hemispheres of the $S^{3}$ bounding the 4-cell.  Note the
distinct rectangles, which cut across the middle of the two 3-balls,
dividing each into two pyramids; one pyramid in each 3-ball has been
made transparent using wire-frame rendering to make the rectangle
visible.  (c) shows a complete partially transparent shaded
rendering of the entire embedded  4-cell projected to 3D, with the outer
octahedron being essentially the equator $S^2$ that is shared by the
two hemispheres (a) and (b) of the $S^3$.

\begin{figure}[!htbp]
 \centerline{
\psfig{figure=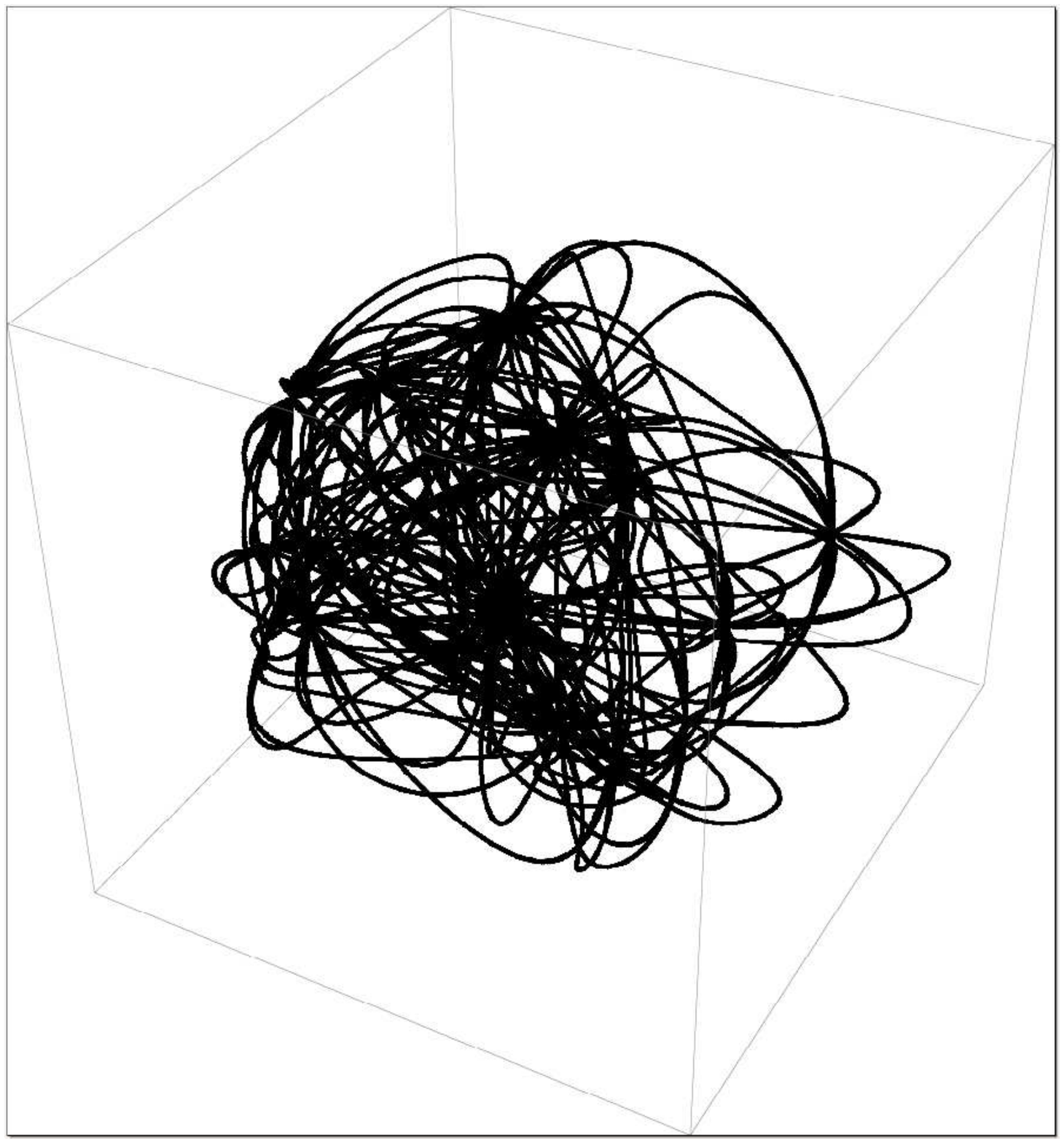,height=18.0cm} }
  \caption[]{}
   \label{K3-edges.fig}
\end{figure}

Finally, we can use computer algebra tools to produce a
representative from each equivalence class of the group $\Gamma_n$
modulo the isotropy group  to generate the explicit forms of all
distinct $12 n^3$ 3-cells and all distinct $6 n^3$ 4-cells. Figure
\ref{K3-edges.fig} represents the K3 surface ($F_{4}$) as a 3D
projection  from its embedding in $\mathbb{R}^{16}$  in terms of all
the edges bounding the $6 n^3 = 384$ 4-cell equators corresponding
to Figure \ref{4cell.fig}(c).  The function of this figure is mainly
to illustrate qualitatively how to use the action of $\Gamma_n$ to
produce the full manifold; more sophisticated interactive
visualization tools are required to expose and explain the
structure, e.g., by interactively selecting and reprojecting subsets
of the tessellation.

 % \bigskip

\end{document}